\documentclass{amsart}
\usepackage{amssymb, amsmath}
\usepackage{bbm, float}
\usepackage{amsfonts}
\usepackage{mathrsfs}
\usepackage{url}
\usepackage{color}

\newtheorem{theorem}{Theorem}[section]
\newtheorem{lemma}[theorem]{Lemma}
\newtheorem{definition}[theorem]{Definition}
\newtheorem{proposition}[theorem]{Proposition}
\newtheorem{example}[theorem]{Example}
\newtheorem{cor}[theorem]{Corollary}
\newtheorem{fact}[theorem]{Fact}
\newtheorem{remark}[theorem]{Remark}
\newtheorem{conjecture}[theorem]{Conjecture}

\def\bref#1{(\ref{#1})}
\def\proof{{\noindent\em Proof:} }
 
\def\Y{{\mathbb{Y}}}
\def\U{{\mathbb{U}}}
\def\L{{\mathbb{L}}}
\def\A{{\mathcal A}}
\def\ml{\mathcal}
\def\mr{\mathscr}
\def\mb{\mathbb}
\def\ff{{\mathcal F}}
\def\ee{{\mathcal E}}

\def\bu{{\mathbf{u}}}
\def\bv{{\mathbf{v}}}
\def\sat{\hbox{\rm{sat}}}
\def\asat{\hbox{\rm{asat}}}
\def\max{\hbox{\rm{max}}}

\def\ord{\hbox{\rm{ord}}}
\def\H{\hbox{\rm{H}}}
\def\lead{\hbox{\rm{ld}}}
\def\sept{\hbox{\rm{S}}}
\def\dim{\hbox{\rm{dim}}}

\def\mod{\hbox{\rm{mod}}}
\def\ord{\hbox{\rm{ord}}}

\def\card{{\rm{card}}}
\def\deg{{\rm deg}}
\def\init{{\rm I}}
\def\rem{{\rm rem}}
\def\rank{{\rm rk}}
\def\trdeg{\hbox{\rm{tr.deg}}}

\def\chow{\rm{Chow}}

\begin{document}

\title[Partial Differential Chow Forms]{Partial Differential Chow Forms and a Type of  Partial Differential Chow varieties}
\author{Wei Li}
\address{KLMM, Academy of Mathematics and Systems Science,
Chinese Academy of Sciences, Beijing 100190, China.}
\email{liwei@mmrc.iss.ac.cn}
\thanks{CONTACT  Wei Li (liwei@mmrc.iss.ac.cn)  KLMM, Academy of Mathematics and Systems Science, Chinese Academy of Sciences, No 55, Zhongguancun East Road, Beijing, 100190 China.}

\begin{abstract}
We first present an intersection theory of irreducible partial differential varieties with quasi-generic differential hypersurfaces.
Then,  we define  partial differential Chow forms for  irreducible partial differential varieties whose Kolchin  polynomials
are of the form $\omega(t)=(d+1){t+m\choose m}-{t+m-s\choose m}$.
And we establish for  partial differential Chow forms most of the basic properties of their ordinary differential counterparts.
Furthermore, we prove that a certain type of partial differential Chow varieties exist.

 \end{abstract}

\maketitle

\section{Introduction}

In their paper on Chow forms \cite{chow}, Chow and van der Waerden described the motivation in these words:
\begin{quote}
It is principally important to represent geometric objects by coordinates. Once this has been done for a specific kind of objects $G$, then it makes sense to speak of an algebraic manifold or an algebraic system of objects $G$,
and to apply the whole theory of algebraic manifolds.
 It is desirable to provide the set of objects $G$ with the structure of an algebraic variety (eventually, after a certain compactification), thus to characterise $G$ by algebraic equations in the coordinates.
\end{quote}
\noindent
Through the theory of Chow forms, they managed to represent projective algebraic varieties or algebraic cycles
by  Chow coordinates;
and Chow further proved that the set of all algebraic cycles of fixed dimension  and degree in the Chow coordinate space is a projective variety, called the Chow variety.

 Chow forms and Chow varieties are basic concepts of algebraic geometry  \cite{chow, hodge}.
 They play an important role in both theoretical and computational aspects of algebraic geometry,
  and  have fruitful applications in elimination theory,  transcendental number theory and algebraic computational complexity theory \cite{eisenbud, Eisenbud-Harris, Complexitychowform, harris1992algebraic,  philippon}.
For instance,  the Chow form  was used by  Brownawell to achieve a major breakthrough in computational algebraic geometry by  proving an effective version of the Hilbert's Nullstellensatz  with optimal bounds    \cite{brownawell}.

Differential algebra, founded by Ritt and Kolchin, is a branch of mathematics analogous to algebraic geometry.
Just as the aim of algebraic geometry is to study solution sets of polynomial equations using algebraic varieties,
the aim of differential algebra is  to study solution sets of algebraic ordinary or partial differential equations using differential varieties \cite{kolchin, ritt}. 
As algebraic equations are special cases of algebraic differential equations, 
algebraic geometry may be viewed as a special case of differential algebra \cite{BC},
 and most of the basic notions of differential algebra are based on those of algebraic geometry. 
Given the importance of the Chow form in algebraic geometry, it is worthwhile to develop a theory of differential Chow forms and differential Chow varieties, and further to see whether they can play similar roles as the algebraic counterparts, for example, whether the bound of the effective differential Nullstellensatz could be improved by applying differential Chow forms.

A systematic development of such a theory was begun by Gao, Li and Yuan in \cite{GLY, LG},
which established a theory of differential Chow forms  for ordinary differential varieties in both affine and projective spaces.
In particular, differential Chow coordinates and two new invariants of ordinary differential varieties (cycles) were introduced.
Take an irreducible differential variety $V$ of differential dimension $d$ for an example. 
Roughly speaking, the differential Chow form of  $V$  is a single differential polynomial $F$
whose general component gives a necessary and sufficient condition when the given $d+1$ differential hyperplanes and $V$ have a common point; and $V$ is determined uniquely by its differential Chow form $F$. 
The coefficient vector of $F$ is called the differential Chow coordinates of $V$.
If the set of all differential cycles of fixed index (dimension, order, and the two new invariants) is a differentially constructible  set in the differential Chow coordinate space,
then we say the differential Chow variety of this index exists.
The existence of ordinary differential Chow varieties was first proved in some special cases by a constructive method in \cite{GLY},
and was finally proved in general cases by Freitag, Li and Scanlon with a model-theoretical proof \cite{FLS}.

However, the theory of   Chow forms has not yet been developed for partial differential varieties.
Unlike the ordinary differential case, an insuperable obstacle is encountered  in the course of defining
 partial differential Chow forms:
due to the more complicated structure of partial differential characteristic sets,
it  is impossible to define differential Chow forms for most of  irreducible partial differential varieties
(see Example \ref{ex-nonexample}).
That is, we may fail to find a single differential polynomial that can represent uniquely the corresponding partial differential variety as in the algebraic and ordinary differential cases.
This leads to the following natural questions: under what conditions can we define partial differential Chow forms?
The set of which kinds of partial differential varieties could be provided  with a structure of partial differential variety (perhaps under Kolchin closure)?
This is what we will deal with in this paper.
Specifically, we will give conditions under which we can define partial differential Chow forms;
under these conditions, we will prove that  partial differential Chow forms have properties similar to those of  their ordinary differential counterparts.
Finally, we will show that partial differential Chow varieties of a certain type  exist.

To define partial differential Chow form, we need to establish a generic intersection theorem for partial differential varieties.
This is also an interesting result in its own right.
The  intersection theorem in algebraic geometry says that
every component of the intersection of two irreducible varieties of dimension $r$ and $s$ in $\mathbb A^n$ has dimension at least $r+s-n$.
However, as pointed out by Ritt, this proposition fails for differential algebraic varieties \cite[p.133]{ritt}.
Recently, Gao, Li and Yuan proved a generic intersection theorem for ordinary differential varieties and generic ordinary differential hypersurfaces \cite{GLY}.
Freitag then generalized this result to the partial differential case using more geometric and model theoretical language \cite{freitag}.
In this paper, we prove the intersection theorem of  partial differential algebraic varieties with quasi-generic partial differential hypersurfaces (to be defined in Definition \ref{def-qgeneric}) using purely differential algebraic arguments.
In particular, when the quasi-generic  differential hypersurface is a generic one, the proof gives more elementary and simplified proofs for  generic intersection theorems either in the ordinary differential case \cite[Theorem 3.6]{GLY} or in the partial differential case \cite[Theorem 3.7]{freitag}.

The  paper is organized as follows.
In section 2, basic notions and preliminary results  are presented.
In section 3, an  intersection theory for quasi-generic partial differential  polynomials  is given.
 In section 4, the definition of the partial differential Chow form and a sufficient condition for its existence are introduced.
Basic properties of the partial differential Chow form are explored  in section 5.
 In section 6, we show partial differential Chow varieties  of a certain type exist.

\section{Preliminaries}

In this section, some basic notation and preliminary results in differential algebra will be given.
For more details about differential algebra, please refer to \cite{BC, kolchin, wu}.

Let $\ff$ be a differential field of characteristic 0 endowed with a finite set of derivations
$\Delta=\{\delta_1,\ldots,\delta_m\}$,
and let $\ml E$ be a fixed universal differential extension field of $\ff$.
 If  $m=1$, $\ff$ and $\ml E$ are called {\em ordinary differential} fields;
 and if $m>1$, they are called {\em partial differential} fields.
Throughout the paper, unless otherwise indicated, all the differential fields (rings) we consider are partial differential fields (rings), and for simplicity, we shall use the prefix ``$\Delta$-" as a synonym of ``partial differential" or ``partial differentially"
when the set of derivation operators in problem are exactly the set $\Delta=\{\delta_1,\ldots,\delta_m\}$.

Let $\Theta$ be the free commutative semigroup (written multiplicatively)  generated by $\delta_1,\ldots,\delta_m.$
Every element $\theta\in\Theta$ is called a {\em derivative operator} and can be expressed uniquely in the form of a product $\prod_{i=1}^m\delta_i^{e_i}$ with $e_i\in\mb N$.
The {\em order} of $\theta$ is defined as $\ord(\theta)=\sum_{i=1}^me_i$.
The identity operator is of order $0$.
For ease of notation, we use $\Theta_s$  to denote the set of all derivative operators of order equal to $s$,
and  $\Theta_{\leq s}$  denotes the set of all derivative operators of order not greater than $s$.
For an element $u\in\ml E$,  denote $u^{[s]}=\{\theta(u): \theta\in\Theta_{\leq s}\}$.

A subset $\Sigma$ of a $\Delta$-extension field $\mathcal {G}$ of
$\mathcal {F}$ is said to be {\em $\Delta$-dependent} over $\mathcal
{F}$ if the set $(\theta\alpha)_{\theta \in \Theta,
\alpha\in\Sigma}$ is algebraically dependent over $\mathcal {F}$,
and is said to be {\em $\Delta$-independent} over $\mathcal {F}$, or
 a family of {\em $\Delta$-$\ff$-indeterminates} in the contrary case.
 In the case $\Sigma$ consists of only one element $\alpha$, we say that
$\alpha$ is $\Delta$-algebraic or $\Delta$-transcendental over
$\mathcal {F}$ respectively.
The {\em $\Delta$-transcendence degree} of $\ml G$ over $\ml F$, denoted by  $\Delta$-$\trdeg\,\ml G/\ml F$,
is the cardinality of any maximal subset of $\ml G$ which are $\Delta$-independent over $\mathcal {F}$.
And the transcendence degree of $\ml G$ over $\ml F$ is denoted by $\trdeg\,\ml G/\ml F$.

Let $\ff\{\Y\}=\ff[\Theta(\Y)]$ be the $\Delta$-polynomial ring over $\ff$ in the  $\Delta$-indeterminates $\Y=\{y_1,\ldots,y_n\}$.
Each element in $\theta(y_i)\in\Theta(\Y)$ is called a derivative, and the order of $\theta(y_i)$ is equal to $\ord(\theta)$.
For a $\Delta$-polynomial $f$ in $\ff\{\Y\}$, the order of $f$ is defined as the maximum of the orders of  all derivatives which appear effectively in $f$, denoted by $\ord(f)$.
A $\Delta$-ideal in $\ff\{\Y\}$ is an ideal which is closed under $\Delta$.
A prime (resp. radical) $\Delta$-ideal is a $\Delta$-ideal which is prime (resp. radical) as an ordinary algebraic ideal.
 Given $S\subset\ff\{\Y\}$, we use $(S)_{\ff\{\Y\}}$ and $[S]_{\ff\{\Y\}}$ to denote the algebraic ideal and the $\Delta$-ideal in $\ff\{\Y\}$ generated by $S$ respectively.

By a $\Delta$-affine space $\mb A^n$, we mean  the set $\ee^n$.
A {$\Delta$-variety} over $\ff$  is  $\mathbb V(\Sigma)=\{\eta\in\ml E^n:\,\, f(\eta)=0, \forall f\in\Sigma\}$ for some set $\Sigma\subseteq\mathcal {F}\{\Y\}$.
 The $\Delta$-varieties in $\mb A^n$ defined over $\ff$ are the closed sets in a topology called the {\em Kolchin topology}.
Given a $\Delta$-variety $V$ defined over $\ff$,  we denote $\mathbb{I}(V)$ to be the set of all $\Delta$-polynomials in
$\ml F\{\Y\}$ that vanish at every point of $V$.  And we have a one-to-one correspondence between $\Delta$-varieties (resp. irreducible $\Delta$-varieties) and radical $\Delta$-ideals (resp. prime $\Delta$-ideal), that is,  for any
$\Delta$-variety $V$  over $\ff$, $\mathbb{V}(\mathbb{I}(V))=V$
and for any radical $\Delta$-ideal $\mathcal{P}$ in $\ff\{\Y\}$,
$\mathbb{I}(\mathbb{V}(\mathcal {P}))=\mathcal {P}$.
For a prime $\Delta$-ideal $\ml P$, a point  $\eta\in \mb V(\ml P)$ is called a {\em generic point} of $\ml P$  (or  $\mb V(\ml P)$) if for any $f\in\ml F\{\Y\}$,  $f(\eta)=0\Longleftrightarrow f\in\ml P$.
A $\Delta$-ideal has a generic point if and only if it is prime.
In this paper, we sometimes use the algebraic version of ideal-variety correspondence, to distinguish from the notation in the differential case,
for an algebraic ideal $\ml P\subseteq\ml F[\Y]$, we use $\ml V(\ml P)$ to denote the algebraic variety in $\mathbb A^n$ defined by $\ml P$;
and for an algebraic variety $V\subseteq\mathbb A^n$, we use $\ml I(V)$ to denote the radical ideal in $\ff[\Y]$ corresponding to $V$.

A homomorphism $\varphi$ from a differential ring $(\mathcal
{R},\Delta)$ to a differential ring $(\mathcal {S}, \Delta')$ with $\Delta'=\{\delta_1',\ldots,\delta_m'\}$ is a
{\em differential homomorphism} if $\varphi\circ \delta_i=\delta_i'\circ
\varphi$ for each $i$. Suppose $\Delta'=\Delta$ and $\mathcal {R}_0$ is a common $\Delta$-subring of
$\mathcal {R}$ and $\mathcal {S}$, $\varphi$ is said to be a $\Delta$-$\mathcal {R}_0$-homomorphism
if  $\varphi$ leaves every element of $\mathcal {R}_0$ invariant.
If, in addition $\mathcal {R}$ is a domain and $\mathcal {S}$ is a $\Delta$-field,
$\varphi$ is called a {\em $\Delta$-specialization} of $\mathcal {R}$ into $\mathcal {S}$.
For $\Delta$-specializations, we have the following lemma which generalizes the similar results both in the ordinary differential case (\cite[Theorem 2.16]{GLY})  and in the algebraic case (\cite[p.168-169]{hodgevol1} and \cite[Lemma 2.13]{GLY}).

\begin{lemma}\label{lem-diffspe}
Let $P_{i} \in \mathcal{F}\{\U, \Y\}\,(i=1, \ldots, \ell)$ be $\Delta$-polynomials in the independent $\Delta$-indeterminates
$\U=(u_{1},  \ldots, u_{r})$ and $\Y$.
Let $\eta$ be an $n$-tuple taken from some extension field of $\mathcal{F}$ free from
$\ff\langle\U\rangle$\footnote{By saying $\eta$ free from
$\ff\langle\U\rangle$, we mean that $\U$ is a set of $\Delta$-$\ff\langle\eta\rangle$-indeterminates.}.
If $P_{i}(\U,\eta)\,(i=1, \ldots, \ell)$ are $\Delta$-dependent over $\ff\langle\U\rangle$,
then for any $\Delta$-specialization $\U$ to $\overline{\U}\in\mathcal{F}^r$, 
$P_{i}(\overline{\U}, \eta)\,(i=1,    \ldots, \ell)$ are $\Delta$-dependent over $\mathcal{F}$.
\end{lemma}
\proof Assume $k=\max_i\ord(P_i)$.
Since $P_{i}(\U,\eta)\,(i=1, \ldots, \ell)$ are $\Delta$-dependent over $\ff\langle\U\rangle$,
there exists $s\in\mb N$ such that the  $\big(P_{i}(\U,\eta)\big)^{[s]}$ are algebraically dependent over $\ff(\U^{[s+k]})$.
When $\U$ $\Delta$-specializes to $\overline{\U}\in\mathcal{F}^r$, $\U^{[s+k]}$ algebraically specializes to $\overline{\U}^{[s+k]}$. 
By \cite[Lemma 2.13]{GLY}, $\big(P_{i}(\overline{\U},\eta)\big)^{[s]}$ are algebraically dependent over $\ff$.
Thus, $P_{i}(\overline{\U}, \eta)\,(i=1,  \ldots, \ell)$ are $\Delta$-dependent over $\mathcal{F}$.
\qed

\subsection{Differential characteristic sets}

A  {\em ranking} of $\ff\{\Y\}$ is a total ordering of the set of derivatives $\Theta(\Y)=\{\theta y_j: j=1,\ldots,n; \theta\in\Theta\}$ 
that satisfies (for any $u,v\in\Theta(\Y)$ and $\delta_k\in\Delta$) the two conditions: 1) $\delta_k u> u$
and 2) $u>v$ $\Rightarrow$ $\delta_k u >\delta_k v$.
Two important kinds of rankings are often used:

    1) {\em Elimination ranking}:
    $y_{i} > y_{j}$ $\Longrightarrow$ $\theta_1  y_{i} >\theta_2 y_{j}$ for any $\theta_1,\theta_2\in\Theta$.

    2) {\em Orderly ranking}: $k>l$  $\Longrightarrow$ for any $\theta_1\in\Theta_k$, $\theta_2\in\Theta_l$ and $i, j$,   $\theta_1  y_{i} >\theta_2  y_{j}$.

   Let $f$ be a $\Delta$-polynomial in $\mathcal {F}\{\Y\}\backslash\ff$ and $\mathscr{R}$  a ranking endowed on it.
   The greatest derivative $\theta y_j$ w.r.t.  $\mathscr{R}$ which  appears effectively in $f$ is called the {\em leader} of $f$, denoted by $\lead({f})$.  Let $d$ be the degree of $f$ in $\lead({f})$. The {\em rank} of $f$ is $\lead({f})^d$, denoted by $\rank(f)$.
   The coefficient of $\rank(f)$ in $f$ is called the {\em initial} of $f$ and denoted by $\init_{f}$.
   The partial derivative of $f$ w.r.t. $\lead({f})$ is called the {\em separant} of $f$, denoted by $\sept_{f}$.
     For any two  $\Delta$-polynomials $f$, $g$ in $\mathcal {F}\{\Y\}\backslash \mathcal {F}$,
    $f$  is said to be of {\em lower rank} than  $g$ if either  $\lead(f)<\lead({g})$ or  $\lead(f)=\lead({g})$ and $\deg(f,\lead(f))<\deg(g,\lead(f))$.
 It is useful to extend the above notion of comparative rank to the whole $\ff\{\Y\}$ by the following convention:
Every element of $\mathcal {F}$ has lower rank than every element of $\mathcal{F}\{\Y\}\backslash\mathcal{F}$
and two elements of $\ff$ have the same rank.

    Let $f$ and $g$ be two $\Delta$-polynomials and $\rank(f)=(\theta y_j)^d$.
    $g$ is said to be {\em partially reduced} w.r.t. $f$ if no proper derivatives of $u_{f}$ appear in $g$.
   $g$ is said to be {\em reduced} w.r.t. $f$ if $g$ is partially reduced w.r.t. $f$ and $\deg(g,\theta y_j)<d$. 
    A set of $\Delta$-polynomials $\mathcal {A}\subset\ff\{\Y\}$ is said to be an {\em autoreduced set} if each $\Delta$-polynomial 
    of $\mathcal {A}$ is reduced  w.r.t.  any other element of $\mathcal{A}$.
   Every autoreduced set is finite.

  Let $\mathcal {A}$ be an autoreduced set. We
    denote $\H_{\mathcal {A}}$ to be  the set of all the initials and
    separants of $\mathcal {A}$ and $\H_{\mathcal {A}}^\infty$ to be the minimal
    multiplicative set containing $\H_{\mathcal {A}}$.
    The {\em $\Delta$-saturation ideal} of $\A$ is defined to be
    $$\sat(\A)=[\mathcal
   {A}]\colon H_{\mathcal {A}}^\infty = \big\{p\in\ff\{\Y\}\big| \exists h\in H_{\mathcal
{A}}^\infty, \,\textrm{ s.t. }\, hp\in[A]\big\}.$$
The algebraic saturation ideal of $\A$ is denoted by $\asat(\A)=(\A)\colon\H_{\mathcal {A}}^\infty$.

    Let $\mathcal {A}=<A_{1},A_{2},\ldots,A_{s}>$ and $\mathcal
    {B}=<B_{1},B_{2},\ldots,B_{l}>$  be two autoreduced
    sets with the $A_{i}$, $B_{j}$ arranged in increasing ordering.
    $\mathcal {A}$ is said to be of {\em lower rank} than $\mathcal {B}$, if
    either 1)\ there is some $k$ ($\leq  \min\{s,l\}$) such that for each
    $i<k$, $A_{i}$ has the same rank as $B_{i}$, and $A_{k}\prec B_{k}$
    or     2)\ $s>l$ and for each $i\in \{1, 2, \ldots, l\}$, $A_{i}$ has the same rank as
    $B_{i}$.  
    If $s=l$ and $A_{i}$ has the same rank as $B_{i}$ for each $i$, we say $\mathcal {A}$ and $\mathcal {B}$ have the same rank.  
    Any sequence of autoreduced sets steadily decreasing in ordering
$\mathcal {A}_{1}\succ \mathcal {A}_{2}\succ \cdots\mathcal
{A}_{k}\succ \cdots$ is necessarily finite.

Let $\mathcal {A}=<A_{1},A_{2},\ldots,A_{t}>$ be an autoreduced
   set with $\sept_{i}$ and $\init_{i}$ as the separant and the initial of $A_{i}$.
Given a $\Delta$-polynomial  $F$, there exists an algorithm, called {\em Ritt's algorithm of reduction}, which reduces
   $F$ w.r.t. $\mathcal {A}$ to a  $\Delta$-polynomial $R$ that is
   reduced w.r.t. $\mathcal {A}$, satisfying the relation
 $$\prod_{i=1}^t\sept_{i}^{d_i}\init_{i}^{e_i} \cdot F \equiv R, \mod \,[\mathcal {A}],$$
   for  $d_{i},e_{i}\in \mb N\, (i=1,2,\ldots,t)$. We call $R$ the {\em remainder} of $F$ w.r.t. $\A$.
   We will need the following result in Section 3.

\begin{proposition}\cite[p.80, Proposition 2]{kolchin} \label{prop-reduction}
Let $\mathcal{A}$ be an autoreduced set of $\ff\{\Y\}$.
If $F_1,\ldots,F_l\in \ff\{\Y\}$, then there exist $\Delta$-polynomials $E_1,\ldots, E_l\in\ff\{\Y\}$,
reduced with respect to $\mathcal{A}$ and of rank no higher than the highest of the ranks of $F_1,\ldots,F_l$,
and there exist $j_A, k_A\in\mathbb N\,(A\in\mathcal{A})$, such that
$$\prod_{A\in\mathcal{A}}\sept_{A}^{j_{A}}\init_{A}^{k_{A}} \cdot F_j \equiv
   E_j, \mod \,[\mathcal {A}]\,\,\,(1\leq j\leq l).$$
\end{proposition}

 Let $\ml J$ be a $\Delta$-ideal in $\ff\{\Y\}$.
 An autoreduced set $\mathcal {C}\subset\ml J$  is said to be a {\em characteristic set} of $\mathcal {J}$,
 if  $\mathcal {J}$ does not contain any nonzero element reduced w.r.t.
$\mathcal {C}$. All the characteristic  sets of $\mathcal {J}$ have the same and minimal rank among all autoreduced sets contained in $\mathcal {J}$. If  $\mathcal {J}$ is prime, $\mathcal {C}$
reduces to zero only the elements of $\mathcal {J}$ and we have
$\mathcal {J}=\sat(\mathcal C)$.
An autoreduced set $\mathcal {C}$ is called {\em coherent} if whenever $A,  A'\in\ml C$ with $\lead(A)=\theta_1(y_j)$ and $\lead(A')=\theta_2(y_j)$ for some $y_j$, the remainder of $\sept_{A'}\frac{\theta}{\theta_1}(A)-\sept_{A}\frac{\theta}{\theta_2}(A')$
w.r.t. $\ml C$ is zero, where $\theta=\text{lcm}(\theta_1,\theta_2)$.\,\,(Here, if $\theta_j=\prod_{i=1}^m\delta_i^{a_{ji}}\,(j=1,2)$ and $\max(a_{1i},a_{2i})=c_i$, then $\theta=\text{lcm}(\theta_1,\theta_2)=\prod_{i=1}^m\delta_i^{c_i}$.)
The following result gives a criterion for an autoreduced set to be a characteristic set of a prime $\Delta$-ideal.

\begin{proposition}\cite[p.167, Lemma 2]{kolchin}
If $\A$ is a characteristic set of a prime  $\Delta$-ideal $\ml P\subset\ff\{\Y\}$, then $\ml P=\sat(\ml A)$, $\ml A$ is coherent, and $\asat(\A)$ is a prime ideal not containing a nonzero element reduced w.r.t. $\ml A$.
Conversely, if $\ml A$ is a coherent autoreduced set of $\ff\{\Y\}$ such that $\asat(\A)$ is a prime ideal not containing a nonzero element reduced w.r.t. $\ml A$, then $\ml A$ is a characteristic set of a prime  $\Delta$-ideal in $\ff\{\Y\}$.
\end{proposition}

\subsection{Kolchin polynomials of prime differential ideals}
Let  $\ml P$ be a prime  $\Delta$-ideal in  $\ff\{\Y\}$ with a generic point $\eta\in\mb A^n$.
The {\em $\Delta$-dimension of $\ml P$}, denoted by $\Delta$-$\dim(\ml P)$, 
is defined as  the $\Delta$-transcendence degree of $\ff\langle\eta\rangle$ over $\ff$. 
A {\em parametric set of $\ml P$}  is a maximal subset $U\subseteq\Y$ such that $\ml P\cap\ff\{U\}=\{0\}$.
Equivalently, the $\Delta$-dimension of $\ml P$ is equal to the cardinality of a  parametric set of $\ml P$.
Let  $\A$ be a characteristic set of $\ml P$ w.r.t. some ranking and
denote $\lead(\A)=\{\lead(F): F\in\A\}$.
 Call $y_j$ a {\em leading variable} of $\A$ if  there exists some $\theta\in\Theta$ such that $\theta(y_j)\in\lead(\A)$;
 otherwise, $y_j$ is called a {\em parametric variable} of $\A$.
 The $\Delta$-dimension of $\ml P$ is also equal to the cardinality of the set of all the parametric variables of $\A$.

For a prime  $\Delta$-ideal, its Kolchin polynomial contains more quantitative information than the $\Delta$-dimension.
To recall the concept of Kolchin polynomial, we need an important numerical polynomial associated to a subset $E\subseteq\mb N^m$.

\begin{lemma} \cite{Kolchin1964, KLMP} \label{lem-numericalpoly}
For every set $E=\{(e_{i1},\ldots,e_{im}): i=1,\ldots,l\}\subseteq \mb N^m$,
let $V_E(t)$ denote the set of all elements $v\in\mb N^m$ such that $v$ is not greater than or equal to any element in $E$ relative to the product order on $\mb N^m$.
Then there exists a univariate numerical polynomial $\omega_E(t)$
such that $\omega_E(t)=\card(V_E(t))$ for all sufficiently large $t.$
Moreover, $\omega_E(t)$ satisfies the following statements:
\begin{itemize}
\item[1)] $\deg(\omega_E)\leq m$, and $\deg(\omega_E)=m$ if and only if $E=\emptyset$. And if $E=\emptyset$, $\omega_E(t)={t+m\choose m}$;
\item[2)] $\omega_E(t)\equiv 0$ if and only if $(0,\ldots,0)\in E$;
\item[3)] If $\min_{i=1}^le_{ik}=0$ for each $k=1,\ldots,m$, then $\deg(\omega_E(t))<m-1$.
\end{itemize}
\end{lemma}

\begin{theorem} \cite[Theorem 2]{Kolchin1964}  \label{th-dimensionpoly}
Let $\ml P$ be a prime $\Delta$-ideal in $\ff\{y_1,\ldots,y_n\}$.
There exists a numerical polynomial $\omega_{\ml P}(t)$ with the following properties:
\begin{enumerate}
\item[1)] For sufficiently large $t\in \mb N$, $\omega_{\ml P}(t)$ equals the dimension of $\ml P\cap\ff[(y_j^{[t]})_{1\leq j\leq n}]$.
\item[2)] $\deg(\omega_{\ml P})\leq m=\card(\Delta)$.
\item[3)] If we write $\omega_{\ml P}(t)=\sum_{i=0}^m a_i{t+i\choose i}$ with $a_i\in \mb Z$, then $a_m$ equals the $\Delta$-dimension of $\ml P$.
\item [4)] If $\ml A$ is a differential characteristic set of $\ml P$ with respect to an orderly ranking on $\ff\{y_1,\ldots,y_n\}$
and if $E_j$ denotes for each $y_j$ the set of points $(l_1,\ldots,l_m)\in\mb N^m$ such that $\delta_1^{l_1}\cdots\delta_m^{l_m}y_j\in\lead(\mathcal A)$, then $\omega_{\ml P}(t)=\sum_{j=1}^n\omega_{E_j}(t)$.
\end{enumerate}
\end{theorem}

The numerical polynomial $\omega_{\ml P}(t)$ is defined to be the {\em Kolchin polynomial} of $\ml P$.
Prime  $\Delta$-ideals  whose characteristic sets consist of a single polynomial are of particular interest to us.
The following result is a partial differential analog of \cite[Lemma 3.10]{GLY} in the ordinary differential case.
 \begin{lemma} \label{le-generalcomponent}
Let $\mathcal{P}$ be a prime $\Delta$-ideal in $\ff\{y_1,\ldots,y_n\}$ and  $A\in\ff\{y_1,\ldots,y_n\}$ an irreducible $\Delta$-polynomial.
Suppose $A$ constitutes a characteristic set of $\mathcal{P}$ under some   ranking $\mathscr{R}$.
Then $\{A\}$ is also a characteristic set of $\mathcal{P}$ under an arbitrary ranking.
In this case, we call $\mathcal{P}$ the {\em general component} of $A$.
\end{lemma}
\proof Follow the proof of \cite[Lemma 3.10]{GLY}.
Suppose $\sept_A$ is the separant of $A$ under $\mathscr{R}$.
Since each element of $\mathcal{P}$ that is partially reduced w.r.t. $A$ is divisible by $A$,
we have $\mathcal{P}=[A]\colon\sept_A^\infty$.
Let  $\mathscr{R}'$ be an arbitrary  ranking and $\theta(y_k)$ be the leader of $A$ under $\mathscr{R}'$.
It suffices to show that there is no nonzero $\Delta$-polynomial in $\mathcal{P}$  reduced with respect to $A$ under $\mathscr{R}'$.
Suppose the contrary and let $f\in\mathcal{P}\backslash\{0\}$  be  reduced with respect to $A$ under $\mathscr{R}'$.
Then $f$ is free from all proper derivatives of $\theta(y_k)$.
Since $f\in \mathcal{P}=[A]\colon\sept_A^\infty$,
 there exist $l\in\mathbb{N}$ and finitely many nonzero polynomials $T_\tau$ for $\tau\in\Theta$ such that
$\sept_A^lf=\sum_{\tau}T_\tau\tau(A)$.
For each $\tau\neq 1$, $\tau(A)=S_A'\cdot\tau\theta(y_k)+L_\tau$, where $S_A'$ is the separant of $A$ under $\mathscr{R}'$.
Substitute $\tau\theta(y_k)=-L_\tau/S_A'$ for each $\tau\neq1$ into both sides of the above identity and remove the denominators, then we get $\sept_A^l \cdot(\sept_{A}')^{l'} f=T_1A$.
Thus, $A$ divides $f$ which implies that $f=0$. 
This contradiction shows that $A$ is also a characteristic set of $\mathcal{P}$ under
any ranking.
\qed
\vskip5pt
Kolchin gave a criterion for a prime $\Delta$-ideal  to be the general component of some $\Delta$-polynomial.
\begin{lemma} \label{lem-iffforgeneralcomponent} \cite[p.160, Proposition 4]{kolchin}
Let $\mathcal{P}\subseteq\ff\{y_1,\ldots,y_n\}$ be a prime $\Delta$-ideal.
Then a necessary and sufficient condition in order that $\mathcal{P}$  be the general component of some polynomial $A$ of order $s$
is 
$$\omega_\mathcal{P}(t)=n{t+m\choose m}-{t+m-s\choose m}.$$
\end{lemma}

The following results on  prime ($\Delta$-)ideals will be used later.

\begin{lemma} \label{lem-dimvariable}
Let $\ml P$ be a prime ideal in the polynomial ring $\ff[x_1,\ldots,x_n]$ of dimension $d>0$.
Assume $\ml P\cap\ff[x_1]=\{0\}$.
Then $\ml J=(\ml P)_{\ff(x_1)[x_2,\ldots,x_n]}$ is a prime ideal of dimension $d-1$.
\end{lemma}
\proof Since $\ml P\cap\ff[x_1]=\{0\}$, $\ml J\neq \ff(x_1)[x_2,\ldots,x_n].$
If $f_1,f_2\in\ff(x_1)[x_2,\ldots,x_n]$ and $f_1f_2\in\ml J$, then there exist $M_1,M_2\in\ff[x_1]$ such that $M_if_i\in\ff[x_1,\ldots,x_n]$ and $M_1f_1M_2f_2\in\ml P$. So either $M_1f_1\in\ml P$ or $M_2f_2\in\ml P$, which implies that either $f_1\in\ml J$
or $f_2\in\ml J$. Thus, $\ml J$ is a prime ideal.

Since $\dim(\ml P)=d$ and $\ml P\cap\ff[x_1]=\{0\}$, without loss of generality, we suppose $\{x_1,x_{2},\ldots,x_{d}\}$ is a parametric set of $\ml P$. We claim that $\{x_{2},\ldots,x_{d}\}$ is a parametric set of $\ml J$.
First, note that $\ml J\cap\ff(x_1)[x_{2},\ldots,x_{d}]=\{0\}$.
For any other variable $x_k\in\{x_{d+1},\ldots,x_n\}$,
$\ml P\cap\ff[x_1,x_{2},\ldots,x_{d},x_k]\neq\{0\}$, so $\ml J\cap\ff(x_1)[x_{2},\ldots,x_{d},x_k]\neq\{0\}$.
Thus, $\{x_{2},\ldots,x_{d}\}$ is a parametric set of $\ml J$, and $\dim(\ml J)=d-1$ follows.
\qed

\begin{lemma} \label{prime-fieldextension}
Let $\ml P$ be a  $\Delta$-prime ideal in $\ff\{\Y\}$ of $\Delta$-dimension $d$.
Suppose $\bu$ is a set of $\Delta$-indeterminates over $\ff$.
Then $[\ml P]_{\ff\langle\bu\rangle\{\Y\}}$ is also a prime $\Delta$-ideal of $\Delta$-dimension $d$.
\end{lemma}
\proof 
Let $\eta$ be a generic point of $\ml P$ free from $\bu$.
It suffices to show that $\eta$ is a generic point of $[\ml P]_{\ff\langle\bu\rangle\{\Y\}}$.
Obviously, $\eta$ is a zero of $[\ml P]_{\ff\langle\bu\rangle\{\Y\}}$.
Suppose $f\in\ff\langle\bu\rangle\{\Y\}$ satisfies $f(\eta)=0$.
By collecting the denominators of $f$, there exists some $D(\bu)\in\ff\{\bu\}$ such that $D(\bu)\cdot f\in\ff\{\bu,\Y\}$.
Write $D(\bu)\cdot f$ in the form
$$D(\bu)\cdot f=\sum_{M}f_M(\Y)M(\bu)$$ 
where $f_M\in\ff\{\Y\}$ and $M(\bu)$'s are distinct $\Delta$-monomials in $\bu$.
Then $f(\eta)=0$ implies that  for each $M$, $f_M(\eta)=0$ and $f_M\in \ml P$.
So $f\in [\ml P]_{\ff\langle\bu\rangle\{\Y\}}$.
Thus, $\eta$ is a generic point of $[\ml P]_{\ff\langle\bu\rangle\{\Y\}}$ and $[\ml P]_{\ff\langle\bu\rangle\{\Y\}}$ is prime.
\qed

\section{Quasi-generic intersection theory for partial  differential polynomials}

In this section, we will prove the quasi-generic intersection theorem with an elementary proof in purely differential algebraic language,
which generalizes generic intersection theorems in both the ordinary differential case \cite{GLY} and  the partial differential case \cite{freitag}.

We recall that a  {\em generic $\Delta$-polynomial in $\mathbb Y=\{y_1,\ldots,y_n\}$ of order $s$ and degree $g$} is a  $\Delta$-polynomial $\mathbb{L}$
of the  form
$$\mathbb{L}=\sum_{M\in\mathscr{M}_{s,g}}u_{M}M(\mathbb Y),$$
where $\mathscr{M}_{s,g}(\mathbb Y)$ is the set of all $\Delta$-monomials  in $\mathbb Y$ of order  $\leq s$ and degree   $\leq g$, 
and all the coefficients $u_{M}\in\mathcal E$ are $\Delta$-$\ff$-indeterminates.  
The $\Delta$-variety $\mathbb V(\mathbb L)\subset\mathbb A^n$ is called a {\em generic $\Delta$-hypersurface}.
If additionally $s=0$ and $g=1$,    $\mathbb V(\mathbb L)$ is called a {\em generic $\Delta$-hyperplane}.

Given $F\in\mathcal F\{\mathbb Y\}$, the set of all $\Delta$-monomials effectively appearing in $F$ is denoted by $\text{supp}(F)$. 
We  now introduce the definition of quasi-generic $\Delta$-polynomials.
\begin{definition} \label{def-qgeneric}
A {\em quasi-generic $\Delta$-polynomial  in $\mathbb Y$ of order $s$}  is a  $\Delta$-polynomial $\mathbb{L}$ of the form
\begin{equation} \label{eq-quasigeneric}
\mathbb{L}=u_0+\sum_{j=1}^nu_jM_j(y_j)+\sum_{M_\alpha\in\text{\rm supp}({\mb L})\backslash\{1,M_1,\ldots,M_n\}}u_{\alpha}M_\alpha(\mathbb Y),
\end{equation}
which satisfies the following conditions:
\begin{itemize}
\item[1)]  for each $j=1,\ldots,n$,  $M_j(y_j)$ is a $\Delta$-monomial in $y_j$ of order $s$;
\vskip5pt
\item[2)] $\{1, M_1(y_1),\ldots,M_n(y_n) \}\subseteq \text{\rm supp}({\mb L})$;
\vskip5pt
\item[3)] the coefficients $u_0, u_j$ and  the $u_\alpha\in\mathcal E$ are $\Delta$-$\ff$-indeterminates.
\end{itemize}
\end{definition}

We give the following main quasi-generic differential  intersection theorem, which generalizes the generic intersection theorem in the ordinary differential case \cite[Theorem 3.6 and Theorem 3.13]{GLY}.
The proof of \cite[Theorem 3.13]{GLY} could not be adapted here due to the more complicated structure of partial differential characteristic sets. 
However, the proof here  could definitely simplify that of its ordinary differential analog.

\begin{theorem} \label{th-quasigeneric-intersection}
Let $V\subseteq\mb A^n$ be an irreducible $\Delta$-variety over $\ff$.
Let $\mathbb{L}$ be a quasi-generic $\Delta$-polynomial of order $s$ with the set of its coefficients $\bu$.
Then
\begin{enumerate}
\item[1)] $V\cap\mathbb V(\mathbb{L})\neq\emptyset$ (over $\ff\langle \bu\rangle$) if and only if  $\Delta$-$\dim(V)>0.$
\vskip5pt
\item[2)] if $\Delta$-$\dim(V)>0$,
then  $V\cap\mathbb V(\mathbb{L})$ is an irreducible $\Delta$-variety over $\ff\langle \bu\rangle$
and its Kolchin  dimension polynomial is
$$\omega_{V\cap \mathbb V(\mathbb{L})}(t)=\omega_V(t)-{t+m-s\choose m}.$$
In particular, the $\Delta$-dimension of $V\cap\mb V(\mb L)$ is equal to $\triangle$-$\dim(V)-1.$
\end{enumerate}
\end{theorem}
\proof
Let $\mathcal{P}=\mathbb{I}(V)\subseteq\ff\{\Y\}$ be the prime $\Delta$-ideal corresponding  to $V$  and $\eta=(\eta_1,\ldots,\eta_n)\in\ee^n$
be a generic point of $\mathcal{P}$  free from $\bu$ (i.e., the $\bu$ are $\Delta$-$\ff\langle\eta\rangle$-indeterminates).
Let $\mb{L}$ be a quasi-generic $\Delta$-polynomial of the form \bref{eq-quasigeneric} and 
set $$\mb T(\Y)=\mb L-u_0=\sum_{j=1}^nu_jM_j(y_j)+\sum_{M_\alpha\in\mr M_{\mb L}\backslash\{1,M_1,\ldots,M_n\}}u_{\alpha}M_\alpha(\Y)\in\ff_1\{\Y\},$$ 
where $\ff_1=\ff\langle\bu\backslash\{u_0\}\rangle$.
Set $\zeta_0=-\mb T(\eta)$.
 
\vskip5pt
1) The proof is similar to that of \cite[Theorem 3.6]{GLY}.
 Let $\ml{J}_0=[\ml P,\mb L]_{\ff_1\{\Y,u_0\}}$.
We first prove that $\ml{J}_0$ is a prime $\Delta$-ideal by showing that $(\eta,\zeta_0)$ is  a generic point of $\ml{J}_0$.
Clearly, $\ml{J}_0$ vanishes at $(\eta,\zeta_0)$.
Given an arbitrary $f\in \ff_1\{\Y,u_0\}$ with $f(\eta,\zeta_0)=0$, we need to show $f\in\ml{J}_0$.
Take the elimination ranking $\mathscr R: y_1<\cdots<y_n<u_0$ of $ \ff_1\{\Y,u_0\}$ and let $f_1$ be the $\Delta$-remainder of $f$ w.r.t. $\mb L$.
Then $f_1\in\ff_1\{\Y\}$ and $f\equiv f_1 \, \mod\, [\mb L]$, which implies that $f_1(\eta)=0$.
By the proof of Lemma \ref{prime-fieldextension}, $\eta$ is also a generic point  of the prime $\Delta$-ideal $[\ml P]_{\ff_1\{\Y\}}$, 
so $f_1\in[\ml P]_{\ff_1\{\Y\}}$ and $f\in \ml{J}_0$ follows. 
Thus,  $\ml{J}_0$ is a prime $\Delta$-ideal with a generic point $(\eta,\zeta_0)$.

Let $\ml J=[\ml P,\mb L]_{\ff\langle\bu\rangle\{\Y\}}$.
We now show that  $\ml J=[1]$ (i.e., $V\cap\mathbb V(\mathbb L)= \emptyset$) if and only if $\Delta$-$\dim(V)=0$.
First suppose $\Delta$-$\dim(V)=0$.
Then for each $j=1,\ldots,n$, $\eta_j$ is $\Delta$-algebraic over $\ff$, and so $\ff_1\langle \eta \rangle$ is $\Delta$-algebraic over $\ff_1$.
Since $\zeta_0\in\ff_1\langle \eta \rangle$, $\zeta_0$ is $\Delta$-algebraic over $\ff_1$.
Thus, $\ml J_0\cap\ff_1\{u_0\}\neq[0]$ and $\ml J=[\ml J_0]_{\ff\langle\bu\rangle\{\Y\}}=[1]$ follows.  
 
 \vskip2pt
For the other direction, suppose $\ml J=[1]$.
 Then $1\in\ml J= [\ml P,\mb L]_{\ff\langle\bu\rangle\{\Y\}}$ implies that $\ml J_0\cap\ff_1\{u_0\}\neq [0]$.
 So $\zeta_0$ is $\Delta$-algebraic over $\ff_1.$
 Consider $-\mb T(\Y)\in\ff\big\{\bu\backslash\{u_0\},\Y\big\}$ and $\zeta_0=-\mb T(\eta)$. 
The quasi-genericness of $\L$ guarantees the existence of a $\Delta$-monomial in $y_j$ for each $j$.
For each $j$, by differentially specializing $u_j$ to $-1$ and all the other elements in $\bu\backslash\{u_0, u_j\}$ to $0$,
$\zeta_0$ will be specialized to $M_j(\eta_j)$.
By Lemma \ref{lem-diffspe}, each $M_j(\eta_j)$, as well as $\eta_j$,  is $\Delta$-algebraic over $\ff$.
So  $\Delta$-$\dim(V)=0$.
Thus,  $\ml J\neq [1]$ if and only if $\Delta$-$\dim(V)>0$.

 \vskip5pt
 2) Assume $\Delta$-$\dim(V)>0$.
 We will show that $\omega_{\ml J}(t)=\omega_V(t)-{t+m-s\choose m}.$
For sufficiently large $t,$ let $\ml I_{t}=\big(\ml P\cap\ff[\Y^{[t]}], \L^{[t-s]}\big)_{\ff_1[\Y^{[t]},u_0^{[t-s]}]}$.
We claim that \begin{itemize}
\item[i)] $\ml I_t\cap\ff_1[u_0^{[t-s]}]=\{0\}$ and $\dim(\ml I_t)=\omega_{\ml P}(t)$;
\vskip5pt
 \item[ii)] $\ml J\cap\ff\langle\bu\rangle[\Y^{[t]}]=(\ml I_{t})_{\ff\langle \bu\rangle[\Y^{[t]}]}$.
\end{itemize}
If i) and ii) are valid, then we have 
\begin{eqnarray}
\omega_{\ml J}(t)&=&\dim(\ml J\cap\ff\langle\bu\rangle[\Y^{[t]}])    \text{ \hskip1.1truecm   (by Definition \ref{th-dimensionpoly})}\nonumber \\
&=& \dim\big((\ml I_{t}\big)_{\ff_1(u_0^{[t-s]})[\Y^{[t]}]})  \text{ \hskip0.8truecm      (by Claim ii))}\nonumber \\
&=& \omega_{\ml P}(t)-{t+m-s\choose m}  \text{\hskip1truecm    (by Lemma \ref{lem-dimvariable} and Claim i)).} \nonumber
\end{eqnarray}
So it remains to show the validity of claims i) and ii).

We first show that $\zeta_t=(\eta^{[t]},\zeta_0^{[t-s]})$ is a generic point of $\ml I_{t}\subset\ff_1[\Y^{[t]},u_0^{[t-s]}]$.
Note that $\ml P\cap\ff[\Y^{[t]}]$ vanishes at $\eta^{[t]}$ and $\L^{(i)}(\zeta_t)=\big(\L(\eta,\zeta_0)\big)^{(i)}=0\,(i\leq t-s)$,
so $\zeta_t$ is a zero of $\ml I_t$.
Suppose $h$ is an arbitrary polynomial in $\ff_1[\Y^{[t]},u_0^{[t-s]}]$ satisfying $h(\zeta_t)=0$.
Let $h_1$ be the pseudo-remainder of $h$ w.r.t. $\L^{[t-s]}$ under the ordering of $\Y^{[t]},u_0^{[t-s]}$ induced by the elimination ranking $\mathscr R$. Then $h_1\in \ff_1[\Y^{[t]}]$ and $h\equiv h_1\,\mod\,(\L^{[t-s]})$.
So $h_1(\eta^{[t]})=0$.
Since $\eta^{[t]}$  is a generic point of the prime ideal $(\ml P\cap\ff[\Y^{[t]}])_{\ff_1[\Y^{[t]}]}$ (by the algebraic version of Lemma \ref{prime-fieldextension}), $h_1\in(\ml P\cap\ff[\Y^{[t]}])_{\ff_1[\Y^{[t]}]}$
and $h\in\ml I_t$ follows.
Thus, $\zeta_t=(\eta^{[t]},\zeta_0^{[t-s]})$ is a generic point of $\ml I_{t}$.

Since $\Delta$-$\dim(V)>0$, by 1), we have $\ml J\neq [1]$ and thus $ \ml J_0\cap\ff_1\{u_0\}=[0]$.
Claim i) follows from the fact that $\ml I_t\cap\ff_1[u_0^{[t-s]}]\subset \ml J_0\cap\ff_1\{u_0\}=[0]$ and 
$\dim(\ml I_t)=\trdeg\,\ff_1(\zeta_t)/\ff_1=\trdeg\,\ff_1(\eta^{[t]})/\ff_1=\omega_{\ml P}(t)$.
Also, note that for any $h\in\ml J\cap \ff_1\{\Y,u_0\}$, there exists $D(u_0)\in\ff_1\{u_0\}$ such that $D\cdot h\in\ml J_0$.
Since $J_0$ is prime and $ \ml J_0\cap\ff_1\{u_0\}=[0]$, $h\in\ml J_0$ and thus we obtain $\ml J\cap \ff_1\{\Y,u_0\}=\ml J_0$.

\vskip4pt
For claim ii), it suffices to show that for each $f\in \ml J\cap\ff\langle\bu\rangle[\Y^{[t]}]$,
$f$ can be written as a  linear combination of polynomials in $\ml P\cap\ff[\Y^{[t]}]$ and $ \L^{[t-s]}$ with coefficients in $\ff\langle\bu\rangle[\Y^{[t]}]$.
 Let $f\in \ml J\cap\ff\langle\bu\rangle[\Y^{[t]}]$.
Multiplying $f$ by some nonzero polynomial in $\ff_1\{u_0\}$ when necessary,
 we can assume $f\in\ff_1[\Y^{[t]},u_0^{[t-s+k]}]$ for some $k\in\mb N$.
 So, $f\in\ml J\cap \ff_1\{\Y,u_0\}=\ml J_0$ and $f(\eta^{[t]},\zeta_0^{[t-s+k]})=0$ follows.
 Let $Z=\cup_{i=1}^{k}\Theta_{t-s+i}$.
 Rewrite $f$ as a polynomial in $\big(\theta(u_0)\big)_{\theta\in Z}$ with coefficients in $\ff_1[\Y^{[t]},u_0^{[t-s]}]$,
  and suppose $$f=\sum_{\alpha} g_\alpha M_\alpha$$ where  $g_\alpha\in\ff_1[\Y^{[t]},u_0^{[t-s]}]$
  and the $M_\alpha$ are finitely many distinct monomials in the variables $\big(\theta(u_0)\big)_{\theta\in Z}$.
So $f(\eta^{[t]},\zeta_0^{[t-s+k]})=0$ implies that $$\sum_{\alpha} g_\alpha(\eta^{[t]},\zeta_0^{[t-s]})M_\alpha\big((\theta(\zeta_0))_{\theta\in Z}\big)=0.$$
If we can show that  $$\big(\theta(\zeta_0)\big)_{\theta\in Z}$$ are algebraically  independent over 
$\ff_1(\eta^{[t]},\zeta_0^{[t-s]})=\ff_1(\eta^{[t]})$,
then  for each $\alpha$, we have $g_\alpha(\eta^{[t]},\zeta_0^{[t-s]})=0$ and $g_\alpha\in  \ml I_t$, which implies that $f\in(\ml I_{t})_{\ff\langle \bu\rangle[\Y^{[t]}]}$.

So it remains to show that $$\big(\theta(\zeta_0)\big)_{\theta\in Z}$$ are algebraically  independent over $\ff_1(\eta^{[t]})$.
Let $\mathcal{A}$ be a $\Delta$-characteristic set of $\ml P$ with respect to some orderly ranking of $\ff\{\Y\}$.
Since $\Delta$-$\dim(V)>0$, there exists at least one $j_0$ such that $y_{j_0}$ is a parametric variable of $\mathcal{A}$.
The quasi-genericness of $\mb L$ guarantees  that $\mb L(\Y)$ effectively involves  a $\Delta$-monomial of order $s$ only in $y_{j_0}$,
say, the term $u_{j_0}M_{j_0}(y_{j_0})$. 
Now consider the polynomials $$\big(\theta(-T(\Y))\big)_{\theta\in Z}\subset\ff\big[(\bu\backslash\{u_0\})^{[t-s+k]},\Y^{[t+k]}\big],$$
and $\big(\theta(-T(\eta))\big)_{\theta\in Z}=\big(\theta(\zeta_0)\big)_{\theta\in Z}.$
Note that  by algebraically specializing $u_{j_0}$ to $-1$ and all the other derivatives of $\bu\backslash\{u_0\}$ to 0, 
$\big(\theta(\zeta_0)\big)_{\theta\in Z}$ are specialized to $\big(\theta(M_{j_0}(\eta_{j_0}))\big)_{\theta\in Z}$. 
If $\big(\theta(\zeta_0)\big)_{\theta\in Z}$ are algebraically dependent over $\ff_1(\eta^{[t]})$,
then by the algebraic version of Lemma \ref{lem-diffspe}, 
$\big(\theta(M_{j_0}(\eta_{j_0}))\big)_{\theta\in Z}$  are algebraically  dependent over $\ff(\eta^{[t]})$.
So there exists a nonzero polynomial $g\big((X_\theta)_{\theta\in Z}\big)\in\ff(\eta^{[t]})[X_\theta:\theta\in Z]$ such that $g\big(\big(\theta(M_{j_0}(\eta_{j_0}))\big)_{\theta\in Z}\big)=0.$ 
By clearing denominators of the coefficients of $g$ (i.e., multiplying some $D(\eta^{[t]})\in\ff[\eta^{[t]}]$) when necessary and replacing $X_\theta$ by $\theta(M_{j_0}(y_{j_0}))$,
we get a nonzero polynomial $$G(\Y)=\sum_{l}g_l(\Y^{[t]})T_l\big(M_{j_0}(y_{j_0})\big)\in\ff\{\Y\}$$ vanishing at $\eta$.
Here, the $T_l(M_{j_0}(y_{j_0}))$ are distinct monomials in $\big(\theta(M_{j_0}(y_{j_0}))\big)_{\theta\in Z}$  
and for each $l$, $g_l\in\ff[\Y^{[t]}]$ does not vanish at $\eta^{[t]}$.
Perform Ritt's algorithm of reduction for all the $g_l$ w.r.t. $\mathcal A$ under the orderly ranking.
By  Proposition \ref{prop-reduction}, there exist  $h_l\in\ff[\Y^{[t]}]\backslash\{0\}$, reduced with respect to $\ml A$, and  natural numbers $j_A, k_A\,(A\in\ml A)$ such that 
$$\prod_{A\in\ml A}\init_A^{j_A}\sept_A^{k_A}\cdot g_l\equiv h_l\,\,\mod\,\,[\ml A],\,\,\,\text{for all $l$'s}.$$
Let $H(\Y)=\sum_{l}h_l(\Y^{[t]})T_l(M_{j_0}(y_{j_0}))\in\ff\{\Y\}$.
Since the order of $M_{j_0}(y_{j_0})$ is  $s$, $\big(\theta(M_{j_0}(y_{j_0}))\big)_{\theta\in Z}$ are algebraically independent over $\ff(\Y^{[t]})$.
Thus, $H(\Y)$ is a nonzero polynomial that is reduced with respect to $\ml A$ and satisfies $H(\eta)=0$,
a contradiction to the fact that $\mathcal A$ is a characteristic set of $\mathcal P$.  Thus, $\big(\theta(\zeta_0)\big)_{\theta\in Z}$ are algebraically  independent over $\ff_1(\eta^{[t]})$ and claim 2) is valid. Consequently, we have proved that $\omega_{V\cap\mathbb{L}}(t)=\omega_{\ml J}(t)=\omega_V(t)-{t+m-s\choose m}.$
\qed

\begin{remark}  
Quasi-generic $\Delta$-polynomials  
are important in that by linear change of coordinates,
``almost all" $\Delta$-polynomials  with a degree-zero term can be transformed to  $\Delta$-polynomials 
with the same supports as  quasi-generic $\Delta$-polynomials. 
Precisely, let $f=a_0+\sum_{M_\alpha\in \text{\rm supp}(f)\backslash\{1\}}a_\alpha M_\alpha(\Y)\in\ff\{\Y\}$ be a $\Delta$-polynomial of order $s$ with $a_0\neq0$. 
Denote $M_\alpha(y_1)\colon=M_\alpha(\Y)\big|_{y_i=y_1,i=1,\ldots,n}$ and set $M_f(y_1)=\max_\alpha\{M_\alpha(y_1)\}$ under  the lexicographical ordering induced by some orderly ranking $\mathscr R$ of $\ff\{\Y\}$. 
Clearly, $\ord(M_f(y_1))=s$.
Let $I=\{\alpha|\, M_\alpha(y_1)=M_f(y_1)\}$.
We associate to $f$ a polynomial $p(x_1,\ldots,x_n)\colon=\sum_{\alpha\in I}a_\alpha\prod_{i=1}^nx_i^{d_{\alpha,i}}\in\ff[x_1,\ldots,x_n]$, where $d_{\alpha,i}=\deg(M_{\alpha},\Theta(y_i))$ for each $i$. 
Assume $p$ is nonzero.
Then under  the linear change of coordinates $\phi\colon y_i=\sum_{k=1}^nb_{ik}z_k\,(i=1,\ldots,n)$ with $b_{ik}\in\ff$  
satisfying $\text{\rm det}(b_{ik})\cdot\prod_{i=1}^np(b_{1i},\ldots,b_{ni})\neq0$,
$\phi(f)\in\ff\{z_1,\ldots,z_n\}$ is a $\Delta$-polynomial of order $s$ with $\{1, M_f(z_1),\ldots,M_f(z_n)\}\subseteq\text{\rm supp}(\phi(f))$ and $\ord(M_f(z_i))=s$.
\end{remark}

{ 
By the proof of Theorem \ref{th-quasigeneric-intersection},
once we know some $y_{i_0}$ which is a parametric variable of a characteristic set of $\mb I(V)$ under some orderly ranking\footnote{Here, the orderly ranking is assumed to guarantee that the order of $h_\ell$ is bounded by $t$ and thus obtain $H(\Y)\neq0$, when following the proof of Theorem \ref{th-quasigeneric-intersection} to prove the corollary.},
 for those $\mb L$ of order $s$ whose support  contains $1$ and a $\Delta$-monomial in $y_{i_0}$ of order $s$ with coefficients $\Delta$-$\ff$-indeterminates, we  can still obtain $\omega_{V\cap\mathbb{L}}(t)=\omega_V(t)-{t+m-s\choose m}$.
  Precisely, we have

\begin{cor}
Let $V\subseteq\mb A^n$ be an irreducible $\Delta$-variety over $\ff$ with $\Delta$-$\dim(V)>0.$
Suppose $\mathcal A$ is a characteristic set of $\mathbb I(V)\subset\ff\{\Y\}$ under some orderly ranking 
with a parametric variable  $y_{i_0}$.
Let $\mathbb{L}(\Y)$ be a $\Delta$-polynomial of order $s$ in the form 
$$\mathbb L=u_0+u_1M_1(y_{i_0})+\sum_{M_\alpha\in\text{\rm supp}({\mb L})\backslash\{1,M_1\}}u_{\alpha}M_{\alpha}(\Y),$$
 where $M_1(y_{i_0})$ is a $\Delta$-monomial in $y_{i_0}$ of order $s$ 
 and the coefficient set $\bu=\{u_0,u_1,u_{\alpha}\}$ is a set of $\Delta$-$\ff$-indeterminates.
Then  $V\cap\mathbb V(\mathbb{L})$ is an irreducible $\Delta$-variety over $\ff\langle \bu\rangle$
and its Kolchin  dimension polynomial is
$$\omega_{V\cap \mathbb V(\mathbb{L})}(t)=\omega_V(t)-{t+m-s\choose m}.$$ 
\end{cor}
}

When $\mb L$ is a generic $\Delta$-polynomial,  Theorem \ref{th-quasigeneric-intersection}  gives the partial differential analog of \cite[Theorem 1.1]{GLY},
which was proven by Freitag with a model-theoretical proof \cite[Theorem 3.7]{freitag}.
\begin{cor} \label{cor-genericintersection}
Let $V$ be an irreducible $\Delta$-variety over $\ff$ with $\omega_V(t)>{t+m\choose m}$.
Let $\mathbb{L}$ be a generic $\Delta$-polynomial of order $s$ and degree $g$ with coefficient set $\bu$.
Then  the intersection of $V$ and $ \mathbb{L}=0$ is a nonempty irreducible $\Delta$-variety over $\ff\langle \bu\rangle$
and its Kolchin   polynomial is
$$\omega_{V\cap\mathbb{L}}(t)=\omega_V(t)-{t+m-s\choose m}.$$
\end{cor}

The following result gives the information of the intersection of several quasi-generic $\Delta$-polynomials, which generalizes \cite[Theorem 3.15]{GLY}.
\begin{cor} \label{cor-dimensionconjecture}
Let $\mb L_i\,(i=1,\ldots,r; r\leq n)$ be independent quasi-generic $\Delta$-polynomials of order $s_i$ respectively.
Suppose $\bu_i$ is the set of coefficients of $\mb L_i$.
Then $V=\mathbb V(\mb L_1,\ldots,\mb L_r)\subset\mathbb A^n$  is an irreducible $\Delta$-variety over ${\ff\langle\bu_1,\ldots,\bu_r\rangle}$ with its Kolchin polynomial equal to
$$\omega_V(t)=(n-r){t+m\choose m}+\sum_{i=1}^r\Big[{t+m\choose m}-{t+m-s_i\choose m}\Big].$$
In particular, if $r=n$, then its  $\Delta$-dimension is 0, the differential type is $m-1$ and the typical $\Delta$-dimension is $\sum_{i=1}^ns_i$.
\end{cor}

\section{Partial Differential Chow forms }
In this section, we  introduce the definition of partial $\Delta$-Chow forms and explore in which conditions on $\Delta$-varieties such that
their $\Delta$-Chow forms exist.

Let $V\subset\mb A^n$ be an irreducible $\Delta$-variety over $\ff$ with $\Delta$-dimension $d$.
Let $$\mathbb L_i=u_{i0}+u_{i1}y_1+\cdots+u_{in}y_n \,\,(i=0,1,\ldots,d)$$
be  independent generic $\Delta$-hyperplanes with coefficients $\bu_i=(u_{i0},\ldots,u_{in})$.
Let 
\begin{equation} \label{eq-J}
\ml J=[\mb I(V), \mb L_0,\ldots,\mb L_d]_{\ff\{\Y,\bu_0,\ldots,\bu_d\}}. 
\end{equation}

\begin{lemma} \label{lemma-J}
$\ml J\cap\ff\{\bu_0,\ldots,\bu_d\}$ is a prime $\Delta$-ideal of codimension 1 
with a parametric set $\bigcup_{i=0}^d \bu_{i}\backslash\{u_{00}\}$.
\end{lemma}
\proof 
Let $\eta=(\eta_1,\ldots,\eta_n)$ be a generic point of $V$ free from each $\bu_i$.
Let $$\zeta_i=-\sum_{k=1}^nu_{ik}\eta_k\,(i=0,\ldots,d) \text{ and } \zeta=(\zeta_0,u_{01},\ldots,u_{0n},\ldots,\zeta_d,u_{d1},\ldots,u_{dn}).$$ 
We first show that $(\eta,\zeta)$ is a generic point of $\ml J$.
Obviously, $\mb I(V)$ and $ \mb L_0,\ldots,\mb L_d$ vanish at $(\eta,\zeta)$.
Suppose $f\in \ff\{\Y,\bu_0,\ldots,\bu_d\}$ with $f(\eta,\zeta)=0$.
Let $f_1$ be the $\Delta$-remainder of $f$ w.r.t. $ \mb L_0,\ldots,\mb L_d$ under some elimination ranking with $\bu<u_{00}<u_{10}<\cdots<u_{d0}$, where  $\bu=\cup_{i=0}^d \bu_{i}\backslash\{u_{i0}\}$. Then $f_1\in\ff\{\Y,\bu\}$ with $f_1(\eta,\zeta)=0$.
Write $f_1=\sum_{M}f_{1M}(\Y)M(\bu)$ as a $\Delta$-polynomial in $\bu$ with coefficients $f_{1M}\in\ff\{\Y\}$, then $f_{1M}(\eta)=0$. So $f_{1M}\in\mb I(V)$ and $f\in \ml J$ follows.
Thus, $(\eta,\zeta)$ is a generic point of $\ml J$ and   $\ml J\cap\ff\{\bu_0,\ldots,\bu_d\}$ is a prime $\Delta$-ideal with a generic point $\zeta$. 

Since the $\Delta$-dimension of $\ml P$ is $d$, 
there exist $d$ of the $\eta_i$ which are $\Delta$-independent over $\ff$. 
Consider $P_{i}=-\sum_{k=1}^nu_{ik}y_k\,(i=1,\ldots, d)$ and the corresponding $\zeta_i$. 
Using Lemma \ref{lem-diffspe} with a contrapositive proof, 
we can prove that $\zeta_1,\ldots,\zeta_d$ are $\Delta$-independent over $\ff\langle\bu\rangle$.
So $\Delta$-tr.deg$\ff\langle\zeta\rangle/\ff\geq (d+1)n+d$.
Note that $\ff\langle\zeta\rangle\subseteq\ff\langle\bu,\eta\rangle$.
Then $\Delta$-tr.deg$\ff\langle\zeta\rangle/\ff=(d+1)n+d$. 
Thus, the codimension of $\ml J\cap\ff\{\bu_0,\ldots,\bu_d\}$ is 1 and 
$\cup_{i=0}^d \bu_{i}\backslash\{u_{00}\}$ is a parametric set of it.
\qed

\vskip3pt
In the ordinary differential case, there always exists a unique irreducible $\delta$-polynomial $F$ such that $\ml J\cap\ff\{\bu_0,\ldots,\bu_d\}$ is  the general component of $F$.
This unique polynomial $F$ is  the $\delta$-Chow form of $V$.
However, unlike the ordinary differential case, for a prime $\Delta$-ideal of codimension 1, it may not be the general component of any single $\Delta$-polynomial, as Example \ref{ex-nonexample}  shows.

\begin{example} \label{ex-nonexample}
Let $m=2$ and $V=\mb V(\delta_1(y),\delta_2(y))\subset\mb A^1$.
Let $\mb L_0=u_{00}+u_{01}y$ and $\ml J=[\mb I(V), \mb L_0]\subset\ml F\{y,u_{00},u_{01}\}$.
Then  $$\ml J\cap\ml F\{u_{00},u_{01}\}={{\sat}\big(u_{01}\delta_1(u_{00})-u_{00}\delta_1(u_{01})},\, u_{01}\delta_2(u_{00})-u_{00}\delta_2(u_{01})\big),$$ which is of codimension 1 but not   the general component of a single $\Delta$-polynomial.
\end{example}

 The above fact makes it impossible to define $\Delta$-Chow forms for all the irreducible
$\Delta$-varieties as in the ordinary differential case \cite[Definition 4.2]{GLY}.  
Below, we define $\Delta$-Chow forms for irreducible $\Delta$-varieties satisfying certain properties.

\begin{definition} \label{def-chowform}
If $\ml J\cap\ff\{\bu_0,\ldots,\bu_d\}$ is the general component of some irreducible  $\Delta$-polynomial $F(\bu_0,\ldots,\bu_d)$, that is,
 $$[\mb I(V), \mb L_0,\ldots,\mb L_d]_{\ff\{\Y,\bu_0,\ldots,\bu_d\}}\cap\ff\{\bu_0,\ldots,\bu_d\}=\sat(F),$$
then we say the  $\Delta$-Chow form of $V$ exists and we call $F$ the {\em  $\Delta$-Chow form} of $V$ or its corresponding
prime $\Delta$-ideal $\mb I(V)$.
\end{definition}

Following this definition, a natural question is to explore in which conditions on $\Delta$-varieties such that their  $\Delta$-Chow forms exist. 
Now, we proceed to give a sufficient condition for the existence of $\Delta$-Chow forms.

\begin{lemma} \label{lm-kolcharset}
Let $\mathcal{P}$ be a prime $\Delta$-ideal in $\ff\{y_1,\ldots,y_n\}$
and $\mathcal A$  a characteristic set of $\ml P$ with respect to some orderly ranking $\mr R$.
Suppose the Kolchin  polynomial of $\ml P$ is
$\omega_{\ml P}(t)=(d+1){t+m\choose m}-{t+m-s\choose m}$ for some $d,s\in\mb N$.
Then $$\lead(\ml A)=\{y_{i_1},\ldots,y_{i_{n-d-1}}, \theta(y_{i_{n-d}})\}$$ for some $\theta\in\Theta_s$ and $n-d$ distinct  variables  $y_{i_1},\ldots,y_{i_{n-d}}$.
\end{lemma}
\proof  For each $j=1,\ldots,n$, let $E_j$ denote the matrix whose row vectors are $(a_1,\ldots,a_m)\in\mb N^m$ such that $\delta_1^{a_1}\cdots\delta_m^{a_m}y_j$ is the leader of an element of $\ml A$. Here, if $y_j$ is not a leading variable,  then set $E_j=\emptyset$.
Suppose the leading variables of $\ml A$ are $y_{i_1},\ldots,y_{i_l}$.
By Theorem \ref{th-dimensionpoly}, $\omega_{\ml P}(t)=\sum_{j=1}^n\omega_{E_j}(t)=(n-l){t+m\choose m}+\sum_{j=1}^l\omega_{E_{i_j}}=(d+1){t+m\choose m}-{t+m-s\choose m}$.    Since $E_{i_j}\neq\emptyset$, the degree of $\omega_{E_{i_j}}$ is less than $m$. Comparing the coefficient of $t^m$ of the both sides of the above equality, we get $l=n-d$.

For $j=1,\ldots,n-d$, let $e_{i_j}=(e_{i_j1},\ldots,e_{i_jn})\in\mb N^m$ be a vector constructed from $E_{i_j}$ with each $e_{i_jk}$  the minimal element of  the $k$-th column of $E_{i_j}$,
and let $H_{i_j}$ be the matrix whose row vectors are the corresponding row vectors of $E_{i_j}$ minus $e_{i_j}$ respectively.
Denote $s_{i_j}=\sum_{k=1}^ne_{i_jk}.$
Then clearly, $\omega_{E_{i_j}(t)}=\omega_{e_{i_j}}(t)+\omega_{H_{i_j}}(t-s_{i_j})$.
By item 3) of Lemma \ref{lem-numericalpoly}, the degree of $\omega_{H_{i_j}}(t-s_{i_j})$ is strictly less than $m-1$.
Thus, $\omega_{\ml P}(t)=(d+1){t+m\choose m}-{t+m-s\choose m}= d{t+m\choose m}+
\sum_{j=1}^{n-d}\omega_{e_{i_j}}(t)+\sum_{j=1}^{n-d}\omega_{H_{i_j}}(t-s_{i_j})$.
So ${t+m\choose m}-{t+m-s\choose m}=\sum_{j=1}^{n-d}[{t+m\choose m}-{t+m-s_{i_j}\choose m}]+\sum_{j=1}^{n-d}\omega_{H_{i_j}}(t-s_{i_j})$. Comparing the coefficients of $t^{m-1}$ and $t^{m-2}$ on the both sides and use the fact ${t+m\choose m}-{t+m-s\choose m}=\frac{s}{(m-1)!}t^{m-1}+\frac{s(m+1)-s^2}{2\cdot(m-2)!}t^{m-2}+o(t^3)$, we get
\[
\left\{ \begin{array}{l}
s=\sum_{j=1}^{n-d}s_{i_j},\\
-s^2/2=-\sum_{j=1}^{n-d}s_{i_j}^2/2+(m-2)!\cdot\sum_{j=1}^{n-d}\text{\rm coeff}\big(\omega_{H_{i_j}},t^{m-2}\big).
\end{array}
\right.
\]
If two of the $s_{i_j}$ are nonzero, then  obviously $-s^2/2<-\sum_{j=1}^{n-d}s_{i_j}^2/2$,
which implies that the above system of equations is not valid.
Thus,  there exists only one $i_j$ such that $s_{i_j}=s$ and all the other $n-d-1$ of the $s_{ij}$ is equal to zero.
Without loss of generality, suppose $s_{i_{n-d}}=s.$
So ${t+m\choose m}-{t+m-s\choose m}= {t+m\choose m}-{t+m-s\choose m}+\sum_{j=1}^{n-d}\omega_{H_{i_j}}(t-s_{i_j})$.
As a consequence, $H_{i_j}=\{(0,\ldots,0)\}$.
Thus, each $E_{i_j}$ has only one row vector, and $\lead(\ml A)=\{y_{i_1},\ldots,y_{i_{n-d-1}}, \theta(y_{i_{n-d}})\}$ for some $\theta\in\Theta_s$.
\qed

\vskip5pt

The following result gives a sufficient condition on $\Delta$-varieties for the existence of $\Delta$-Chow forms.

\begin{theorem} \label{th-chowsufficient}
Let $V\subseteq\mb A^n$ be an irreducible $\Delta$-variety over $\ff$ with  $$\omega_V(t)=(d+1){t+m\choose m}-{t+m-s\choose m}$$ for some $s\in\mb N$.
Then the $\Delta$-Chow form of $V$ exists. And the order of the $\Delta$-Chow form of $V$ is $s$.
\end{theorem}

\proof
Let $\mathcal{P}=\mb I (V)\subset \ff\{\Y\}$.
Let $\mathcal P^{\star}=[\ml P, \mb L_1,\ldots,\mb L_d]_{\ff\langle \bu_1,\ldots,\bu_d\rangle\{\Y\}}$.
Then by Corollary \ref{cor-genericintersection},
$\ml P^{\star}$ is a prime $\Delta$-ideal of $\Delta$-dimension 0 and $\omega_{\ml P^\star}={t+m\choose m}-{t+m-s\choose m}$.
Let $\mathcal J_0=[\ml P^\star, \mb L_0]_{\ff\langle\bu_1,\ldots,\bu_d\rangle\{\Y,\bu_0\}}$.
Recall that $\ml J=[\ml P, \mb L_0,\ldots,\mb L_d]_{\ff\{\Y,\bu_0,\ldots,\bu_d\}}$ as introduced in (\ref{eq-J}).
Clearly, $\ml J_0=[\ml J]_{\ff\langle\bu_1,\ldots,\bu_d\rangle\{\Y,\bu_0\}}$.
Let $$\ml I_0=\ml J_0\cap \ff\langle\bu_1,\ldots,\bu_d\rangle\{\bu_0\} \text{  and  } \ml I=\ml J\cap\ff\{\bu_0,\ldots,\bu_d\}.$$ 
By Lemma \ref{lemma-J},  both $\ml I_0$ and $\ml I$ are prime $\Delta$-ideals of codimension 1 with parametric sets $\bu_0\backslash\{u_{00}\}$ and  $\bigcup_{i=0}^d\bu_i\backslash\{u_{00}\}$ respectively.
For any $f\in\ml I_0\cap \ff\{\bu_0,\ldots,\bu_d\}$, there exists $M\in\ff\{\bu_1,\ldots,\bu_d\}$ such that $M f\in\ml I$.
Since $M\notin \ml I$, $f\in\ml I$.
Thus, $\ml I=\ml I_0\cap\ff\{\bu_0,\ldots,\bu_d\}$. 
We claim that
($\ast$) $\ml I$ is the general component of some irreducible $F\in\ff\{\bu_0,\ldots,\bu_d\}$  if and only if $\ml I_0$ is the general component of some irreducible $G\in\ff\langle\bu_1,\ldots,\bu_d\rangle\{\bu_0\}$. 

To show the validity of claim ($\ast$),  first suppose $\ml I$ is the general component of some irreducible $F\in\ff\{\bu_0,\ldots,\bu_d\}$.
Fix two elimination rankings for $\ff\{\bu_0,\ldots,\bu_d\}$ and $\ff\langle\bu_1,\ldots,\bu_d\rangle\{\bu_0\}$ 
with $u_{00}$ higher than any other variable.
Since $\bigcup_{i=0}^d\bu_i\backslash\{u_{00}\}$ is a parametric set of $\ml I$,
$\lead(F)=\theta (u_{00})$ for some $\theta$.  
Given any $g\in \ml I_0$, let  $g_1$ be the partial remainder of $g$ w.r.t. $F$.
Then $g_1\in\mathcal I_0.$
If $g_1\neq0$, there exists $L\in\ff\{\bu_1,\ldots,\bu_d\}$ such that $Lg\in \ml I$.
Since $F$ is a characteristic set of $\ml I$ and $Lg$ is partially reduced w.r.t. $F$,
$Lg$ is divisible by $F$ over $\ff$, and thus $g$ is divisible by $F$ over $\ff\langle\bu_1,\ldots,\bu_d\rangle$.
Thus, $F$ constitutes a characteristic set of $\ml I_0$ and $\ml I_0$ is the general component of $F$.
For the other direction, suppose $\ml I_0$ is the general component of some irreducible $G\in\ff\langle\bu_1,\ldots,\bu_d\rangle\{\bu_0\}$. 
Clearly, $\lead(G)=\theta(u_{00})$ for some $\theta$ and 
there exists $N\in\ff\{\bu_1,\ldots,\bu_d\}$ such that $NG\in\ff\{\bu_0,\bu_1,\ldots,\bu_d\}$ is irreducible.
Since $\ml I_0=[G]\colon\left(\frac{\partial G}{\partial \theta(u_{00})}\right)^{\infty}$, we have $\ml I=\ml I_0\cap\ff\{\bu_0,\ldots,\bu_d\}=\left([NG]\colon\big(N\frac{\partial G}{\partial \theta(u_{00})}\big)^{\infty}\right)_{\ff\{\bu_0,\ldots,\bu_d\}}=\sat(NG)$.
So $\ml I$ is the general component of $NG$ and  ($\ast$) is proved.

Thus, by ($\ast$),  it suffices to consider for the case $\dim(V)=0$, that is, to show the $\Delta$-Chow form of $V$ exists if $\omega_V(t)={t+m\choose m}-{t+m-s\choose m}$ for some $s\in\mb N$.
Now suppose $\dim(V)=0$ and let $\eta=(\eta_1,\ldots,\eta_n)$ be a generic point of $V$ free from $\bu_0$.
Let $\zeta_0=-\sum_{j=1}^nu_{0j}\eta_j$.
Then by the proof of Lemma \ref{lemma-J}, $(\zeta_0,u_{01},\ldots,u_{0n})$ is a generic point of $\ml I=[\mb I(V), \mb L_0]\cap\ff\{\bu_0\}$.

On the one hand, since $\zeta_0^{[t]}\subseteq\ff(u_{01}^{[t]},\ldots,u_{0n}^{[t]},\eta^{[t]})$,  we have
$$\omega_{(\zeta_0, u_{01},\ldots,u_{0n})}(t)\leq \omega_{(u_{01},\ldots,u_{0n},\eta)}(t)=(n+1){t+m\choose m}-{t+m-s\choose m}.$$
On the other hand, by Lemma \ref{lm-kolcharset}, $\omega_V(t)={t+m\choose m}-{t+m-s\choose m}$  implies that
the leading variables of a characteristic set of $\ml A$ with respect to an orderly ranking is $\{y_{i_1},\ldots,y_{i_{n-1}},\theta(y_{i_n})\}$ with $\theta\in\Theta_s$. So $\{\tau(\eta_{i_n}):\,\tau\in\Theta_{\leq t},\theta\nmid \tau\}$ is algebraically independent over $\ff$.
By  the contrapositive of the algebraic version of Lemma \ref{lem-diffspe}, $S:=\{\tau(\zeta_{0}):\,\tau\in\Theta_{\leq t},\theta\nmid\tau\}$ is algebraically independent over $\ff(u_{01}^{[t]},\ldots,u_{0n}^{[t]})$. Note that $\card(S)={t+m\choose m}-{t+m-s\choose m}.$ Thus,  we have
\begin{eqnarray}
&&\omega_{(\zeta_0, u_{01},\ldots,u_{0n})}(t)=\trdeg\,\ff(u_{01}^{[t]},\ldots,u_{0n}^{[t]},\zeta_0^{[t]})/\ff\nonumber\\ &=&\trdeg\,\ff(u_{01}^{[t]},\ldots,u_{0n}^{[t]})/\ff
+\trdeg\,\ff(u_{01}^{[t]},\ldots,u_{0n}^{[t]})(\zeta_0^{[t]})/\ff(u_{01}^{[t]},\ldots,u_{0n}^{[t]})\nonumber\\ &\geq& n{t+m\choose m}+{t+m\choose m}-{t+m-s\choose m}. \nonumber
\end{eqnarray}

Thus, $\omega_{(\zeta_0, u_{01},\ldots,u_{0n})}(t)=(n+1){t+m\choose m}-{t+m-s\choose m}.$
By Lemma \ref{lem-iffforgeneralcomponent}, there exists an irreducible $\Delta$-polynomial $F$ of order $s$ such that $\ml J=\sat(F)$, so the
$\Delta$-Chow form of $V$ exists and it is of order $s$.
\qed

\vskip5pt

We conjecture that for the existence of $\Delta$-Chow form of $V$, the Kolchin polynomial $\omega_V(t)=(d+1){t+m\choose m}-{t+m-s\choose m}$  is   also a  necessary condition:
\begin{conjecture} \label{conjecture}
Let $V\subseteq\mb A^n$ be an irreducible $\Delta$-variety  over $\ff$ of differential dimension $d$.
Then  the $\Delta$-Chow form of $V$ exists if and only if
$$\omega_V(t)=(d+1){t+m\choose m}-{t+m-s\choose m}\,\,\text{ for some $s\in\mb N$}.$$
\end{conjecture}

\vskip5pt
In the remaining sections of the paper, we focus on irreducible $\Delta$-varieties $V\subseteq\mathbb A^n$ of Kolchin  polynomial $\omega_V(t)=(d+1){t+m\choose m}-{t+m-s\choose m}$ for some $s\in\mathbb N$ whose $\Delta$-Chow forms exist guaranteed by Theorem \ref{th-chowsufficient}.

\vskip3pt
Given $d+1$ vectors $\mathbf{a}_i=(a_{i0},a_{i1},\ldots,a_{in})\in\mathbb A^{n+1}(\ee)$ for $i=0,\ldots,d$, 
let $$\mathbb L_i(\mathbf{a}_i)=a_{i0}+a_{i1}y_1+\cdots+a_{in}y_n$$ be the corresponding hyperplanes over $\ff\langle \mathbf{a}_i\rangle$. 
A natural question in  (differential) intersection theory arises as follows: 
\begin{center}
Under which conditions can we have $ V\cap\mathbb L_0(\mathbf{a}_0)\cap\cdots\cap\mathbb L_d(\mathbf{a}_d)\neq\emptyset$?
\end{center}
In the algebraic case, Chow and van der Waerden showed that the vanishing of the Chow form gives a necessary and sufficient condition such that the projective variety $V$ and $d+1$ projective hyperplanes have a common point \cite{chow}.
 Similarly, we also have a geometric interpretation for partial differential Chow forms.
 The following result shows that the general component of the differential Chow form of $V$  gives a necessary and sufficient condition in the Kolchin closure sense such that $V$ and the given $d+1$ $\Delta$-hyperplanes have a nonempty intersection.
\begin{proposition}
Let $V\subset\mb A^n$ be an irreducible $\Delta$-variety  of Kolchin   polynomial
$\omega_V(t)=(d+1){t+m\choose m}-{t+m-s\choose m}$ and $F(\bu_0,\ldots,\bu_d)$ be the $\Delta$-Chow form of $V$.
Let $$S=\left\{(\mathbf{a}_0,\ldots,\mathbf{a}_d)\in (\mathbb A^{n+1})^{d+1}\,\,\big| \,\,V\cap\mathbb L_0(\mathbf{a}_0)\cap\cdots\cap\mathbb L_d(\mathbf{a}_d)\neq\emptyset\right\}.$$
Then the Kolchin closure of $S$ is the general component of $F$.
\end{proposition}
\proof Let $W$ be the Kolchin closure of $S$ in $(\mathbb A^{n+1})^{d+1}$.
Recall that 
$$\mathcal J=[\mathbb I(V),\mathbb L_0,\ldots,\mathbb L_d]_{\ff\{\Y,\bu_0,\ldots,\bu_d\}} \text{  and  }  \mathcal J\cap\ff\{\bu_0,\ldots,\bu_d\}=\sat(F).$$
For any $(\mathbf{a}_0,\ldots,\mathbf{a}_d)\in S$, there exists $\xi\in\mathbb A^n$ such that $\xi\in V\cap\mathbb L_0(\mathbf{a}_0)\cap\cdots\cap\mathbb L_d(\mathbf{a}_d)$.
Then $(\xi,\mathbf{a}_0,\ldots,\mathbf{a}_d)\in\mathbb V(\mathcal J)$ and so $(\mathbf{a}_0,\ldots,\mathbf{a}_d)\in \mathbb V(\sat(F))$.
Thus, $S\subseteq\mathbb V(\sat(F))$ and $W\subseteq \mathbb V(\sat(F))$.
On the other hand, by the proof of Lemma \ref{lemma-J},  $(\eta,\zeta)$ is a generic point of $\mathcal J$, so $\zeta$ is a generic point of $\mathbb V(\sat(F))$.
Clearly, $\eta\in V\cap\mathbb L_0(\zeta_{u,0})\cap\cdots\cap\mathbb L_d(\zeta_{u,d})$, where $\zeta_{u,i}=(\zeta_i,u_{i1},\ldots,u_{in})$.
So $\zeta\in S\subset W$. 
Thus, $\mathbb V(\sat(F))\subseteq W$.
Hence, we have $\mathbb V(\sat(F))=W$.
\qed

\vskip5pt
Below is an example of $\Delta$-Chow forms.

\begin{example} Let $\Delta=\{\delta_1,\delta_2\}$.
Let $\ml P=[\delta_1(y_1),y_2-y_1^2]\subset\ff\{y_1,y_2\}$.
Clearly, $\omega_{\ml P}(t)={t+2\choose 2}-{t+1\choose 2}=t+1$.
The $\Delta$-Chow form of $\ml P$ is
$$F(\bu_0)=\delta_1(u_{00})^2u_{02}^2-2\delta_1(u_{00})u_{00}\delta_1(u_{02})u_{02}-\delta_1(u_{00})\delta_1(u_{01})u_{01}u_{02}+u_{00}^2(\delta_1(u_{02}))^2+ $$$$\delta_1(u_{00})u_{01}^2\delta_1(u_{02})
+ u_{00}(\delta_1(u_{01}))^2u_{02}-u_{00}\delta_1(u_{01})u_{01}\delta_1(u_{02})+
 \delta_1(u_{00})\delta_1(u_{01})u_{01}u_{02}.$$

\end{example}

\section{Properties of the partial differential Chow form}
In this section, we will prove basic properties of $\Delta$-Chow forms.
In particular, we will show the $\Delta$-Chow form is $\Delta$-homogenous and  has  a Poisson-type product formula similar to its ordinary differential counterpart.

\subsection{Partial differential Chow forms are differentially homogenous}

In this section, we will show that the $\Delta$-Chow form is  $\Delta$-homogenous.
Recall that $\ff$ is a $\Delta$-field with the set of derivations $\Delta=\{\delta_1,\ldots,\delta_m\}$
and  $\Theta$ is the set of all derivative operators.
Given two derivatives $\theta_1=\prod_{i=1}^m\delta_i^{a_i}$ and $\theta_2=\prod_{i=1}^m\delta_i^{b_i}\in\Theta$,
if $a_i\leq b_i$ for each $i$,  then we denote $\theta_1|\theta_2$.
In case $\theta_1|\theta_2$,
  we denote $\frac{\theta_2}{\theta_1}=\prod_{i=1}^m\delta_i^{b_i-a_i}$, and denote the product of binomial coefficients $\prod_{i=1}^m{b_i\choose a_i}$ by ${\theta_2\choose \theta_1}$. It is easy to verify that $\theta(fg)=\sum_{\tau|\theta}{\theta\choose \tau}\cdot\frac{\theta}{\tau}(f)\cdot\tau(g)$ for all $f,g\in\ff$.

\begin{definition} \label{def-homogenous}
A $\Delta$-polynomial $f\in\ff\{y_0,y_1,\ldots,y_n\}$ is said to be {\em $\Delta$-homogenous} of degree $r$ if $f(\lambda y_0,\lambda y_1,\ldots,\lambda y_n)=\lambda^rf(y_0,y_1,\ldots,y_n)$ holds for a $\Delta$-indeterminate $\lambda$ over $\ff\{y_0,y_1,\ldots,y_n\}.$
\end{definition}

The following lemma is a partial differential analog of the Euler's criterion on homogenous polynomials, which was listed as an exercise
in \cite[p.71]{kolchin}.
\begin{proposition} \label{lem-diffhom}
A necessary and sufficient condition that $f\in\ff\{y_0,y_1,\ldots,y_n\}$ be $\Delta$-homogenous of degree $r$
 is that $f$ satisfies the following system of equations:
 \begin{equation} \label{eq-homegenous}\sum_{\tau\in\Theta}\sum_{j=0}^n{\tau\theta\choose \theta}\tau(y_j)\cdot\frac{\partial f}{\partial \tau\theta(y_j)}=
  \left\{
  \begin{array}{lll}
   rf,&\quad&\theta=1 \\
   0,&\quad&\theta\in\Theta,\,\theta\neq1.
 \end{array}
 \right.
 \end{equation}
\end{proposition}

\proof
Denote $\Y=(y_0,\ldots,y_n)$ temporarily for convenience.
Let $\lambda$ be a $\Delta$-indeterminate over $\ff\{\Y\}$.

First, we show the necessity.
Suppose $f$ is $\Delta$-homogenous of degree $r$.
Then $f(\lambda \Y)=\lambda^rf(\Y)$.
Differentiating both sides of this equality w.r.t. $\theta(\lambda)$, we get
\begin{eqnarray}
&&\sum_{\tau\in\Theta}\sum_{j=0}^n{\tau\theta\choose \theta}\tau(y_j)\cdot\frac{\partial f}{\partial \tau\theta(y_j)}(\lambda \Y)=\sum_{j=0}^n \sum_{\tau\in\Theta}\frac{\partial \tau\theta(\lambda y_j)}{\theta(\lambda)}\frac{\partial f}{\partial \tau\theta(y_j)}(\lambda \Y)
\nonumber\\&=&\frac{\partial f(\lambda\Y)}{\partial \theta(\lambda)}=\left\{
\begin{array}{lll}
rf(\Y)\lambda^{r-1},&&\theta=1\\
0,&&\theta\in\Theta, \theta\neq1.
\end{array}\right.\nonumber\end{eqnarray}
Setting $\lambda=1$,  we obtain \bref{eq-homegenous}.

To  show the sufficiency, suppose \bref{eq-homegenous} holds.
We will show $f(\lambda \Y)=\lambda^rf(\Y)$ for some $r$.
Let $\theta\in\Theta$ satisfy ($\ast$) $\frac{\partial}{\partial \tau\theta(\lambda)}f(\lambda\Y)=0$ for all $\tau\in\Theta\backslash\{1\}$.
Then we have
\begin{eqnarray} \label{eq-sufficiency}
\lambda\cdot\frac{\partial}{\partial \theta(\lambda)}f(\lambda\Y)&=&\sum_{\tau\in\Theta}{\tau\theta\choose \theta}\tau(\lambda)\frac{\partial}{\partial \tau\theta(\lambda)}f(\lambda\Y)\\
&=&\sum_{\tau\in\Theta}{\tau\theta\choose \theta}\tau(\lambda)\sum_{j=0}^n\sum_{\xi\in\Theta}{\xi\tau\theta\choose \tau\theta}\xi(y_j)\frac{\partial}{\partial \xi\tau\theta(y_j)}f(\lambda\Y) \nonumber \\
&=& \sum_{j=0}^n\sum_{\tau\in\Theta}\sum_{\xi\in\Theta,\xi |\tau}{\frac{\tau}{\xi}\theta\choose \theta}\frac{\tau}{\xi}(\lambda){\tau\theta\choose \frac{\tau}{\xi}\theta}\xi(y_j)\frac{\partial}{\partial \tau\theta(y_j)}f(\lambda\Y) \nonumber \\
 \nonumber \\
 &=& \sum_{j=0}^n\sum_{\tau\in\Theta}{\tau\theta\choose \theta}\Big(\sum_{\xi\in\Theta,\xi |\tau}{\tau\choose \xi}\frac{\tau}{\xi}(\lambda)\xi(y_j)\Big)\frac{\partial}{\partial \tau\theta(y_j)}f(\lambda\Y) \nonumber \\
 \nonumber \\
&=&\sum_{j=0}^n\sum_{\tau\in\Theta}{\tau\theta\choose \theta}\tau(\lambda y_j)\frac{\partial}{\partial \tau\theta(y_j)}f(\lambda\Y) \nonumber\\
&=& \Big(\sum_{\tau\in\Theta}\sum_{j=0}^n{\tau\theta\choose \theta}\tau(y_j)\cdot\frac{\partial f}{\partial \tau\theta(y_j)}\Big)_{\Y=\lambda\Y}. \nonumber 
\end{eqnarray}
For each $\theta\neq1$ satisfying ($\ast$), by \bref{eq-homegenous} and \bref{eq-sufficiency},  
we obtain $\frac{\partial}{\partial \theta(\lambda)}f(\lambda\Y)=0$.
Since $\frac{\partial}{\partial \tau(\lambda)}f(\lambda\Y)=0$ for all $\tau\in\Theta_{>\ord(f)}$, each $\theta\in\Theta_{\ord(f)}$ satisfies ($\ast$).
By backward induction on $\Theta(\lambda)$, 
we have $\frac{\partial}{\partial \theta(\lambda)}f(\lambda\Y)=0$ for each $\theta\in\Theta\backslash\{1\}$.
Thus, we can take $\theta=1$ and get $\lambda\cdot\frac{\partial}{\partial  \lambda}f(\lambda\Y)=rf(\lambda\Y)$ from \bref{eq-homegenous} and \bref{eq-sufficiency}.
So $$\frac{\partial \lambda^{-r}f(\lambda\Y)}{\partial \lambda}=-r \lambda^{-r-1}f(\lambda\Y)+ \lambda^{-r}\frac{\partial f(\lambda\Y)}{\partial \lambda}=0,$$ and $f(\lambda\Y)=\lambda^rf(\Y)$ follows. Thus, $f(\Y)$ is $\Delta$-homogenous of degree $r$.
\qed

\vskip5pt
 Now, we show that the $\Delta$-Chow form $F(\bu_0,\ldots,\bu_d)$ is  $\Delta$-homogeneous of the same degree in  each $\bu_i$,
 which is a partial differential analog of \cite[Theorem 4.17]{GLY}.
\begin{theorem} \label{th-diffhomogenous}
Let $V\subset\mb A^n$ be an irreducible $\Delta$-variety of Kolchin  polynomial
$\omega_V(t)=(d+1){t+m\choose m}-{t+m-s\choose m}$.
Let $F(\bu_0,\ldots,\bu_d)$ be the $\Delta$-Chow form of $V$.
Then $F(\bu_0,\ldots,\bu_d)$ is $\Delta$-homogenous of the same degree $r$  in  each $\bu_i$.
\end{theorem}
\proof  By the definition of $\Delta$-Chow form, $F(\bu_0,\ldots,\bu_d)$ has the symmetric property in the sense that
interchanging $\bu_i$ and $\bu_j$ in $F$, the resulting polynomial and $F$ differ at most by a sign.
In particular, $F$ is of the same degree in each $\bu_i$.
So it suffices to show the $\Delta$-homogeneity of $F$ for $\bu_0$.

Let $\eta=(\eta_1,\ldots,\eta_n)$ be a generic point of $V$. 
For  $i=0,\ldots,d$, let $$\zeta_i=-\sum_{j=1}^nu_{ij}\eta_j  \text{\,\,  and    }\,\, \zeta_{u,i}=(\zeta_i,u_{i1},\ldots,u_{in}).$$ 
By the proof of Lemma \ref{lemma-J},
$(\zeta_{u,0},\zeta_{u,1},\ldots,\zeta_{u,d})$ is a generic point of 
$$[\mathbb I(V),\mb L_0,\ldots,\mb L_d]_{\ff\{\Y,\bu_0,\ldots,\bu_d\}}\cap\ff\{\bu_0,\ldots,\bu_d\}=\sat(F).$$ 
Here, a  ranking $\mathscr R$ of $\ff\{\bu_0,\ldots,\bu_d\}$ is assumed.
Let $\lambda$ be a $\Delta$-indeterminate over $\ff\langle\bu_0,\ldots,\bu_d,\Y\rangle$.
Similarly as in the proof of Lemma \ref{lemma-J}, we can also show that 
$(\lambda\zeta_{u,0},\zeta_{u,1},\ldots,\zeta_{u,d})$ is also a generic point of $\sat(F)$.
So $F(\lambda\zeta_{u,0},\zeta_{u,1},\ldots,\zeta_{u,d})=0$.
Let $F(\lambda\bu_0,\bu_1,\ldots,\bu_d)=\sum_{M}F_MM(\lambda)$ 
where the $M(\lambda)$'s are distinct $\Delta$-monomial in $\lambda$
and $F_M\in\ff\{\bu_0,\ldots,\bu_d\}$.
Since $$F(\lambda\zeta_{u,0},\zeta_{u,1},\ldots,\zeta_{u,d})=\sum_{M}F_M(\zeta_{u,0},\zeta_{u,1},\ldots,\zeta_{u,d})M(\lambda)=0,$$
we have $F_M(\zeta_{u,0},\zeta_{u,1},\ldots,\zeta_{u,d})=0$ for each $M$.
So $F_M\in\sat(F).$ 
Since $F_M$ is partially reduced w.r.t. $F$ and $\deg(F_M)\leq F$,
there exists $a_M\in\ff$ such that $F_M=a_MF$.
Thus, $F(\lambda\bu_0,\bu_1,\ldots,\bu_d)=\left(\sum_{M}a_MM(\lambda)\right)\cdot F$.
Regard the two $\Delta$-polynomials at both sides of this equality as polynomials in $\ff\langle\lambda\rangle\{\bu_0,\ldots,\bu_d\}$,
by comparing their coefficients of the leading monomial under the lexicographic ordering induced by $\mathscr R$,
we obtain $\sum_{M}a_MM(\lambda)=\lambda^r$ for some $r\in\mathbb N$.
Thus, $F$ is $\Delta$-homogenous in $\bu_0$ of degree $r$.
 \qed

\begin{definition}
The number $r$ in Theorem \ref{th-diffhomogenous} is defined to be the {\em $\Delta$-degree} of the  $\Delta$-variety $V$
or its corresponding prime $\Delta$-ideal.
\end{definition}

\subsection{Factorization of partial differential Chow forms}
In this section, we follow the techniques in \cite[Sect. 4.4]{GLY} to derive  Poisson-type product formulae for partial differential Chow forms.
For this purpose, fix an orderly ranking $\mr R$ on $\bu_0,\ldots,\bu_n$ with $u_{00}$ greater than any other $u_{ij}$.
 Suppose $\lead(F)=\theta(u_{00})$ and $\theta$ is reserved for this derivative operator temporarily in this section.
 Let $$\ff_\bu=\ff\langle\bu_1,\ldots,\bu_d,u_{01},\ldots,u_{0n}\rangle\,\, \text{ and }\,  \ff_0=\ff_\bu\big(\tau(u_{00}):\tau\in\Theta, \theta\nmid\tau\big).$$
 Regard $F$ as a univariate  polynomial $f\big(\theta(u_{00})\big)$ in  $ \theta(u_{00})$ with coefficients in $\ff_0$
 and suppose $g=\deg(F,\theta(u_{00}))$.
Then $f\big(\theta(u_{00})\big)$ is irreducible over $\ff_0$ and in a suitable algebraic extension field of $\ff_0$,
$f(\theta(u_{00}))=0$ has $g$ roots $\gamma_{1},\ldots,\gamma_{g}$.
Thus
\begin{equation}\label{eq-fac00}
f(\theta(u_{00}))=A(\bu_{0},\bu_{1},\ldots,\bu_{d})\prod^g_{l=1}\big(\theta(u_{00})-\gamma_{l}\big)
\end{equation}
where  $A(\bu_{0},\bu_{1},\ldots,\bu_{d})\in\mathcal
{F}\{\bu_0,\ldots,\bu_d\}$ is free from $\theta(u_{00})$.

 For each $l=1, \ldots, g$, let
 \begin{equation}\label{eq-ftau}
 \mathcal{F}_l=\mathcal {F}_0\big(\gamma_{l}\big)\end{equation}
 be an algebraic extension of $\mathcal {F}_0$ defined by
 $f\big(\theta(u_{00})\big)=0$.
We will define derivations $\delta_{l,1},\ldots,\delta_{l,m}$ on $\mathcal{F}_l$ so
that $\big(\mathcal{F}_l ,\{\delta_{l,1},\ldots,\delta_{l,m}\}\big)$ becomes a partial differential field.
This can be done step by step in a very natural way.
For the ease of notation,  for each  $\tau=\prod_{k=1}^m\delta_{k}^{d_k}$ with $(d_1,\ldots,d_m)\in\mb N^m$, we denote $\tau_l=\prod_{k=1}^m\delta_{l,k}^{d_k}$.
In step 1, for each $a\in \mathcal{F}_\bu$, define $\tau_l(a)=\tau(a)$, in particular, $\delta_{l, k}(a) = \delta_k(a)$ for each $k=1,\ldots,m$.
In step 2, we need to define the derivatives of $u_{00}$.
For all $\tau\in\Theta$ with $\theta\nmid \tau$ or $\tau=\theta$,  define $\tau_l(u_{00})$ as follows:
\[
\tau_l(u_{00})=\left\{\begin{array}{lll}
\tau(u_{00})\in\ff_0(\subseteq\ff_l),&& \theta\nmid\tau\\
\gamma_l\in\ff_l,&& \tau=\theta.\\
\end{array}
\right.
\]
And for all $\tau\in\Theta$ with $\theta|\tau$ and $\tau\neq\theta$,
we define $\tau_l(u_{00})$ inductively on the ordering of $\Theta(u_{00})$ induced by $\mathscr R$.
Since $F$, regarded as a univariate polynomial $f$ in $\theta(u_{00})$, is a minimal polynomial of $\gamma_{l}$,
$\sept_{f}=\frac{\partial f}{\partial \theta(u_{00})}$ does not vanish at $\theta(u_{00})=\gamma_{l}$.
First, for the minimal $\tau=\delta_k\theta$ for some $k\in\{1,\ldots,m\}$, define
$$\tau_l(u_{00})=\delta_{l,k}(\gamma_l)\stackrel{\mathrm{def}}{=}-T/\sept_{f}\big|_{\theta(u_{00})=\gamma_l},$$
where $\delta_k(f)=S_f\cdot\delta_k\theta(u_{00})+T$.
 This is reasonable, since all the derivatives of $u_{00}$ involved in $S_f$ and $T$ have been defined in the former steps and we should have  $\delta_{l,k}\big(f(\gamma_l)\big)=\sept_{f}\big|_{\theta(u_{00})=\gamma_l}\delta_{l,k}(\gamma_l)+T\big|_{\theta(u_{00})=\gamma_l}=0$.
 Suppose all the derivatives of $u_{00}$ less than $\tau(u_{00})=\prod_{k=1}^m\delta_{k}^{d_k}\theta(u_{00})$ have been defined, we can proceed in the similar way to define $\tau_l(u_{00})=\prod_{k=1}^m\delta_{l,k}^{d_k}(\gamma_l)$.
 Namely, use the differential polynomial $\tau(f)=S_f\cdot\tau(u_{00})+T_{\tau}$ and define
 $\tau_l(u_{00})=-T_\tau/\sept_{f}\big|_{\pi\theta(u_{00})=\pi(\gamma_l), \,\pi\theta<\tau}$.
In this way,  $(\ff_l,\{\delta_{l,1},\ldots,\delta_{l,m}\})$ is a partial differential field which
can be considered as a finitely differential extension field of $(\ff_{\bu},\Delta)$.

Since $\ff_{\bu}$ is a finitely generated $\Delta$-extension field of $\ff$ contained in $\ee$.
By the definition of universal differential extension fields, there exists a $\Delta$-extension field
$\ff^*\subset\ee$ of $\ff_{\bu}$ and a differential $\ff_{\bu}$-isomorphism $\varphi_l$ from $(\ff_l,\{\delta_{l,1},\ldots,\delta_{l,m}\})$ to
$(\ff^*,\Delta)$. For a polynomial $G\in\ff\{\Y\}$ and a point $\eta\in\ff_l^n$,
$G(\eta)=0$ implies $G(\varphi_l(\eta))=0$.  For convenience, by saying
$\eta$ is in a $\Delta$-variety $V$ over $\ff$, we mean $\varphi_l(\eta)\in V$.
Summing up the above results, we have a partial differential analog of \cite[Lemma 4.24]{GLY}.

\begin{lemma}\label{lm-ftau}
$(\ff_l,\{\delta_{l,1},\ldots,\delta_{l,m}\})$ is a finitely differential
extension field of  $(\ff_{\bu}, \Delta)$, which is
differentially $\ff_{\bu}$-isomorphic to a differential subfield of
$\ee$. \end{lemma}

Note that  the above defining steps give a differential homomorphism $\phi_l$ from $(\ff\{\bu_0,\ldots,\bu_d\}, \Delta)$ to the differential field $(\ff_l,\{\delta_{l,1},\ldots,\delta_{l,m}\})$ for each $l$ by mapping $\tau(u_{ij})$ to $\tau_l(u_{ij})$.
That is, for  a $\Delta$-polynomial $p\in\mathcal{F}\{\bu_0, \ldots,\bu_d\}$,
$\phi_l(p)$ is  obtained from $p$ by substituting $\tau\theta(u_{00})=\tau_l(\gamma_l)$.
Then we have the following result similar to \cite[Lemma 4.25]{GLY}.

\begin{lemma}\label{lm-ftaup}
Let $P\in\mathcal{F}\{\bu_0, \ldots,\bu_d\}$.
Then $P\in\sat(F)$ if and only if $\phi_l(P)=0$.
\end{lemma}
\proof If $P\in\sat(F)$, then there exists $m\in\mathbb N$ such that $S_F^mP\in[F]$.
Since $\phi_l$ is a differential homomorphism and $\phi_l(F)=0$, $\phi_l(S_F^mP)=0$.
Note from the above that $\phi_l(S_F)\neq0$, so $\phi_l(P)=0$ follows.
For the other side, suppose $\phi_l(P)=0$.
Let $R$  be the differential remainder of $P$ w.r.t. $F$ under the ranking $\mr R$.
Since $\phi_l(P)=0$, $\phi_l(R)=0$.
Note that $R$ is free from all the proper derivatives of $\theta(u_{00})$ and $\deg(R,\theta(u_{00}))<g$.
So $R|_{\theta(u_{00})=\gamma_l}=0$, which implies from the irreducibility of $F$ that $R$ is divisible by $F$.
Thus, $R=0$ and $P\in\sat(F)$.
\qed

\begin{remark}\label{rm-ftau}
Similar to the ordinary differential case, in order to make $\ff_l$ a partial differential field, we need to
introduce derivations $\delta_{l,1},\ldots,\delta_{l,m}$  related to $\gamma_l$ and
there does not exist a unique set of  derivations  to make all $\ff_l(l=1,\ldots,g)$ differential fields.
 \end{remark}

Below, we give the following Poisson-type product formula, which is a partial differential analog of \cite[Theorem 4.27]{GLY}.
\begin{theorem}\label{th-fac1}
Let  $F(\bu_{0},\bu_{1},\ldots,\bu_{d})$ be the $\Delta$-Chow form of an irreducible $\Delta$-variety over $\ff$ of
Kolchin  polynomial $\omega_V(t)=(d+1){t+m\choose m}-{t+m-s\choose m}$.
Fix an orderly ranking with $u_{00}>u_{ij}$ and suppose $\lead(F)=\theta(u_{00})$  and $g=\deg\big(F,\theta(u_{00})\big)$.
Then, there exist $\xi_{l
1},\ldots,\xi_{l n}$ in a differential extension field
$(\ff_l,\{\delta_{l,1},\ldots,\delta_{l,m}\})$  of $(\mathcal
{F}_{\bu},\Delta)$ such that
\begin{eqnarray}\label{eq-fac10}
F(\bu_{0},\bu_{1},\ldots,\bu_{d})=A(\bu_{0},\bu_{1},\ldots,\bu_{d})\prod^g_{l=1}\theta\left(u_{00}+ \sum_{\rho=1}^n
u_{0\rho}\xi_{l \rho}\right)
\end{eqnarray}
where  $A(\bu_{0},\bu_{1},\ldots,\bu_{d})$ is in $\mathcal
{F}\{\bu_0,\ldots,\bu_d\}$.
Note that equation \bref{eq-fac10} is formal and should be
understood in the following precise meaning:
$\theta(u_{00}+\sum_{\rho=1}^n u_{0\rho}\xi_{l \rho})
\stackrel{\triangle}{=}\theta(u_{00})+\theta_l(\sum_{\rho=1}^n
u_{0\rho}\xi_{l \rho}\big).$
\end{theorem}
\proof Follow the proof of \cite[Theorem 4.27]{GLY} and the  notation above.
By Lemma \ref{lm-ftaup}, $\phi_l(S_F)\neq0$.
Let $\xi_{lj}=\phi_l(\frac{\partial F}{\partial \theta(u_{0j})})\big/\phi_l(S_F)$ for $j=1,\ldots,n$ and $\xi_l=(\xi_{l1},\ldots,\xi_{ln})\in\ff_l^n$.
We will prove \[\gamma_{l}=-\theta_l\Big(\sum_{j=1}^nu_{0j}\xi_{lj}\Big).\]

 Let $\xi=(\xi_1,\ldots,\xi_n)$ be a generic point of $V$ and $\zeta_i=-\sum_{j=1}^nu_{ij}\xi_j$.
 Then $F(\zeta_0,u_{01},\ldots,u_{0n};\ldots;\zeta_d,u_{d1},\ldots,u_{dn})=0$.
Differentiating the equailty  w.r.t.
$\theta(u_{0j})$ on both sides, we have
\begin{equation} \label{eq-partial}\overline{\frac{\partial F}{\partial \theta(u_{0j})}}+
 \overline{\frac{\partial F}{\partial \theta(u_{00})}}\cdot(-\xi_{j})=0,\end{equation}
 where the $\overline{\frac{\partial F}{\partial \theta(u_{0j})}}$ are
obtained by substituting $\zeta_{i}$ to $u_{i0}\,(i=0, \ldots,
d)$ in $\frac{\partial F}{\partial \theta(u_{0j})}$.
 Multiplying $u_{0j}$ to the above equation and for $j$ from 1 to $n$, adding them together, we have
\[\sum_{j=1}^n u_{0j}\overline{\frac{\partial F}{\partial \theta(u_{0j})}}+
 \overline{\frac{\partial F}{\partial \theta(u_{00})}}\cdot(-\sum_{j=1}^n u_{0j}\xi_{j})
 =\sum_{j=1}^n u_{0j}\overline{\frac{\partial F}{\partial \theta(u_{0j})}}+
 \overline{\frac{\partial F}{\partial \theta(u_{00})}}\cdot\zeta_0=0.\]
Thus,  $\sum_{j=0}^n u_{0j}\frac{\partial F}{\partial
\theta(u_{0j})} \in\,\sat(F)$.
By Lemma \ref{lm-ftaup},  $$\sum_{j=1}^n u_{0j}\phi_l\big(\frac{\partial F}{\partial
\theta(u_{0j})} \big)+\phi_l(u_{00})\phi_l\big(\frac{\partial F}{\partial
\theta(u_{00})} \big)=0,$$ so $\phi_l(u_{00})=-\sum_{j=1}^nu_{0j}\xi_{lj}$.
Thus, $\theta_l(\phi_l(u_{00}))=\phi_l(\theta(u_{00}))=\gamma_l=-\theta_l(\sum_{j=1}^nu_{0j}\xi_{lj})$.
 Substituting them into
equation~\bref{eq-fac00}, \bref{eq-fac10} is proved.\qed

\vskip5pt
We have the following interesting result  similar to \cite[Theorem 4.34]{GLY}.
\begin{theorem} \label{th-zerochowform}
The points $(\xi_{l1},\ldots,\xi_{l n})\,(l=1,\ldots,g)$ in \bref{eq-fac10} are  generic points of the $\Delta$-variety $V$ over $\ff$.
If $d>0$, they also satisfy the equations \[u_{\sigma0}+\sum_{\rho=1}^n u_{\sigma \rho}y_{\rho}=0\,(\sigma=1,\ldots,d).\]
\end{theorem}
\proof  Follow the proof of  \cite[Theorem 4.34]{GLY}.
Suppose $P(y_{1},\ldots,y_{n}) \in \mathcal {F}\{\Y\}$ is
any $\Delta$-polynomial vanishing on $V$. Then
$P(\xi_{1},\ldots,\xi_{n})=0$. From (\ref{eq-partial}),
$\xi_{\rho}=\overline{\frac{\partial f}{\partial
 \theta(u_{0\rho})}}\big/\overline{\frac{\partial{f}}{\partial\theta(u_{00})}}$,
so we have \[P\left(\overline{\frac{\partial F}{\partial
\theta(u_{01})}}\bigg/\overline{\frac{\partial{F}}{\partial
\theta(u_{00})}}, \ldots,\overline{\frac{\partial F}{\partial
\theta(u_{0n})}}\bigg/\overline{\frac{\partial{F}}{\partial
\theta(u_{00})}}\right)=0,\]
where $\overline{\frac{\partial F}{\partial \theta(u_{0\rho})}}$ are
obtained by substituting $\zeta_{i}$ to $u_{i0}\,(i=0, 1, \ldots,
d)$ in $\frac{\partial f}{\partial  \theta(u_{0\rho})}$.
Thus, there exists an $t \in \mathbb{N}$, such that
$$\big(\frac{\partial{F}}{\partial
\theta(u_{00})}\big)^t\cdot P\Big( {\frac{\partial F}{\partial
\theta(u_{01})}}\bigg/ {\frac{\partial{F}}{\partial
\theta(u_{00})}}, \ldots, {\frac{\partial F}{\partial
\theta(u_{0n})}}\bigg/ {\frac{\partial{F}}{\partial
\theta(u_{00})}}\Big)\in\sat(F).$$
By Lemma \ref{lm-ftaup}, we have
$P(\xi_{l1},\ldots,\xi_{l n})=0$,  which means that
$(\xi_{l1},\ldots,\xi_{l n}) \in V$.

Conversely, for any $Q \in \mathcal {F}\{\Y\}$ such that
$Q(\xi_{l1},\ldots,\xi_{l n})=0$, by Lemma \ref{lm-ftaup}, there exists $t \in
\mathbb{N}$ s.t. $\widetilde{Q}=(\frac{\partial F}{\partial
\theta(u_{00})})^t Q(\frac{\partial F}{\partial
\theta(u_{01})}\Big/\frac{\partial F}{\partial
\theta(u_{00})},\ldots,\frac{\partial F}{\partial
\theta(u_{0n})}\Big/\frac{\partial f}{\partial \theta(u_{00})}) \in \sat(F)$.
So $Q(\xi_{1},\ldots,$ $\xi_{n})=0$. Thus, $(\xi_{l1},\ldots,\xi_{l n})$ is a generic point of $V$.

By equation (\ref{eq-partial}),  $\overline{\frac{\partial F}{\partial \theta(u_{0j})}}+
 \overline{\frac{\partial F}{\partial \theta(u_{00})}}\cdot(-\xi_{j})=0$,
so we have $\sum_{j=1}^n u_{\sigma j}\overline{\frac{\partial
F}{\partial
\theta(u_{0j})}}+\zeta_{\sigma}\overline{\frac{\partial{F}}{\partial \theta({u}_{00})}}=0$.
Thus, $\sum_{j=0}^n u_{\sigma j}\frac{\partial F}{\partial
\theta(u_{0j})}\in\sat(F)$. If $\sigma \neq 0$, then $\sum_{j=0}^n u_{\sigma j}\phi_l(\frac{\partial F}{\partial
\theta(u_{0j})})$ $=0$. Consequently,
$u_{\sigma0}+\sum_{j=1}^n u_{\sigma j}\xi_{l j}=0\,(\sigma=1,\ldots,d)$. \qed

\begin{remark}
In \cite{GLY}, the number $g$ in the Poission-type formula is defined as the leading differential degree of   $\delta$-cycles,
which has similar geometric meaning as the degree of algebraic varieties. 
But  in the partial differential case, the leading differential degree could not be defined, for the number $g$ in Theorem
\ref{th-fac1} depends on the orderly ranking we choose to get the Poisson-type product formula.
 {For example, let  $\Delta=\{\delta_1,\delta_2\}$ and $V\subset\mathbb A^1$ be the general component of $A=\delta_1(y_1)(\delta_2(y_1))^2+1\in\ff\{y_1\}$.
 Then the differential Chow form of $V$ is $F(\bu_0)=-u_{01}^3\delta_1(u_{00})(\delta_2(u_{00}))^2+2u_{00}u_{01}^2\delta_1(u_{00})\delta_2(u_{00})\delta_2(u_{01})-u_{00}^2\delta_1(u_{00})u_{01}(\delta_2(u_{01}))^2+
 u_{00}u_{01}^2\delta_1(u_{01})(\delta_2(u_{00}))^2-2u_{00}^2u_{01}\delta_1(u_{01})\delta_2(u_{00})\delta_2(u_{01})
 +u_{00}^3\delta_1(u_{01})(\delta_2(u_{01}))^2+u_{01}^6$.
 There are two orderly rankings of $\ff\{y_1\}$, that is, $\mathscr R_1: \delta_1^{i_1}\delta_2^{i_2}> \delta_1^{j_1}\delta_2^{j_2}\Leftrightarrow (i_1+i_2,i_1)>_{lex}(j_1+j_2,j_1)$ and $\mathscr R_2: \delta_1^{i_1}\delta_2^{i_2}> \delta_1^{j_1}\delta_2^{j_2}\Leftrightarrow (i_1+i_2,i_2)>_{lex}(j_1+j_2,j_2)$.
 Under $\mathscr R_1$, $\text{ld}(F)=\delta_1(u_{00})$ and $g=1$;
 while under $\mathscr R_2$, $\text{ld}(F)=\delta_2(u_{00})$ and $g=2$.
 }
Also, even under a fixed orderly ranking, 
the leaders of the $\Delta$-Chow forms of two irreducible $\Delta$-varieties with the same Kolchin polynomial may be distinct, 
so it is difficult to define partial differential algebraic cycles through $\Delta$-Chow forms as we did in the ordinary differential case.
\end{remark}

\vskip2pt

We conclude this section by proposing the following properties of  $\Delta$-Chow forms, 
which are similar to \cite[Lemma 4.10]{GLY} and \cite[Lemma 3.9]{FLS} in the ordinary differential case
and will be used in section 6.
\begin{theorem} \label{th-relationchow}
Let $V\subseteq\mb A^n$ be an irreducible $\Delta$-variety of Kolchin polynomial $\omega_V(t)=(d+1){t+m\choose m}-{t+m-s\choose m}$ and $F(\bu_0,\ldots,\bu_d)$  the $\Delta$-Chow form of $V$. The following assertions hold.
\begin{itemize}
\item[1)] Let $\mathscr R$ be some elimination ranking   satisfying $u_{ij}<u_{00}<y_1<\cdots<y_n$.
 Let $\lead(F)=\theta(u_{00})$ and $\sept_F$ the separant of $F$.
 Then $$\{F,\sept_Fy_1-\frac{\partial F}{\partial \theta(u_{01})},\ldots,\sept_Fy_n-\frac{\partial F}{\partial \theta(u_{0n})}\}$$ is a characteristic set of $[\mb I(V), \mb L_0,\ldots,\mb L_d]_{\ff\{\Y, \bu_0,\ldots,\bu_d\}}$ w.r.t. $\mathscr R$.

\item[2)]
Given $(\bv_0,\ldots,\bv_d)\in(\mathbb A^{n+1})^{d+1}$,
if $F(\bv_0,\ldots,\bv_d)=0$ and $\sept_F(\bv_0,\ldots,\bv_d)\neq0$,
then $V$ and $v_{i0}+v_{i1}y_1+\cdots+v_{in}y_n=0\,(i=0,\ldots,d)$ have at least one point in common.

\item[3)] $\left(\mb I(V)\cap\ff[\Y^{[s]}], \mb L_0^{[s]}, \ldots, \mb L_d^{[s]}\right)\cap\ff[\bu_0^{[s]},\ldots,\bu_d^{[s]}]=\big(F(\bu_0,\ldots,\bu_d)\big).$
\end{itemize}
\end{theorem}
\proof The proof of item 1) is similar to  \cite[Lemma 4.10]{GLY}. And item 2) is a direct consequence of item 1).

 3) By Theorem \ref{th-chowsufficient}, $\ord(F)=s$, so $(F)=[\mb I(V), \mb L_0,\ldots,\mb L_d]\cap\ff[\bu_0^{[s]},\ldots,\bu_d^{[s]}].$
By the proof of Lemma \ref{lemma-J}, $\zeta^{[s]}$ is a generic point of $(F)$.
Similarly, we can show
$(\eta^{[s]},\zeta^{[s]})$ is a generic point of $\ml J_s:=\big(\mb I(V)\cap\ff[\Y^{[s]}], \mb L_0^{[s]}, \ldots, \mb L_d^{[s]}\big)$ in $\ff[\Y^{[s]},\bu_0^{[s]},\ldots,$ $\bu_d^{[s]}]$.
Thus, $\zeta^{[s]}$ is also a generic point of $\ml J_s\cap\ff[\bu_0^{[s]},\ldots,\bu_d^{[s]}]$ and 3) follows.
\qed

\section{Partial differential Chow varieties of a certain type exist}

As mentioned in the introduction, to study a specific kind of geometric objects,
it is important and useful to represent them by coordinates and further show that the set of objects is actually an algebraic system.
For us, this specific kind of objects are irreducible $\Delta$-varieties with Kolchin  polynomial $(d+1){t+m\choose m}-{t+m-s\choose m}$.
 As in the ordinary differential case,  we could represent these $\Delta$-varieties by coordinates.
 \begin{definition}
Let $V$ be an irreducible $\Delta$-variety over $\ff$ of Kolchin polynomial $\omega_V(t)=(d+1){t+m\choose m}-{t+m-s\choose m}$ and of $\Delta$-degree $r$.
Let $F(\bu_0,\ldots,\bu_d)$ be the $\Delta$-Chow form of $V$.
The coefficient vector of $F$, regarded as a point in a higher dimensional projective space determined by $n$ and $(d, s, r)$,   is defined to be the {\em $\Delta$-Chow coordinate} of $V$.
 \end{definition}

  Fix $n$ and an index $(d,s,r)$.   Let $G_{(n,d,s,r)}$ be a functor from the category of $\Delta$-fields of characteristic $0$ to the category of sets,
 which associates each $\Delta$-field $\ml F$ of characteristic 0 with
 the set $G_{(n,d,s,r)}(\ml F)$, consisting of all irreducible $\Delta$-varieties $V\subseteq\mb A^n$  over  $\ml F$ with $\omega_V(t)=(d+1){t+m\choose m}-{t+m-s\choose m}$ and $\Delta$-degree $r$.

 \begin{definition}
  If  $G_{(n,d,s,r)}$ is represented by some $\Delta$-constructible set over $\mathbb Q$,
  meaning that there is a $\Delta$-constructible set defined over $\mb Q$ 
  and a natural isomorphism between the functor $G_{(n,d,s,r)}$ and the functor given by
 this $\Delta$-constructible set  (regarded also as a  functor from the category of $\Delta$-fields of characteristic $0$ to the category of sets),
 then we call this $\Delta$-constructible set the
 \emph{$\Delta$-Chow variety} of index $(d, s, r)$ of $\mb A^n$, and denote it by $\Delta$-Chow$(n,d,s,r)$.
 In this case, we also say that the $\Delta$-Chow variety of $\mb A^n$ of  index $(d,s,r)$  exists.
 \end{definition}

In this section, we will show that $\Delta$-Chow varieties $\Delta$-Chow$(n,d,s,r)$  exist for all chosen $n, d, s, r$.
Similar to the ordinary differential case, 
the main idea is to first definably embed $G_{(n,d,s,r)}$ into a finite disjoint union $\ml C$ of the chosen algebraic Chow varieties
and then show the image of $G_{(n,d,s,r)}$ is a definable subset of $\ml C$.
So the language from model theory of partial differentially closed fields (see \cite{mcgrail, pierce, sanchez}) will be used.
Assume $\ml E$ is a saturated $\Delta$-closed field of characteristic $0$ (i.e., $\ml E \models  \text{DCF}_{0,m}$) and $\mb A^n=\ml E^n$ throughout this section.

\subsection{Definable properties and Prolongation admissible varieties}

Here are some basic notions and results from model theory to be used in the proof of the main theorem.
For more details and explantations, see \cite{FLS}.

We say that a
family of sets $\{ X_{a} \}_{a \in B}$ is a \emph{definable family} if there
are formulae $\psi (x; y)$ and $\phi(y)$ so that $B$ is the set of realizations
of $\phi$ and for each $a \in B$,  $X_{a}$ is the set of realizations
of $\psi(x;a)$.

Given a property $\mathrm P$ of definable sets, 
we say that $\mathrm P$ is definable in families if for any family of definable sets $\{ X_a \}_{a \in B}$
given by the formulae $\psi(x;y)$ and $\theta(y)$, 
there is a formula $\phi(y)$ so that the set
$\{ a \in B ~:~ X_a \text{ has property } \mathrm P \}$ is defined by $\phi$.

Given an operation $\mathcal L$ which takes a set and returns another set, we
say that $\ml L$ is definable in families if for any
family of definable sets $\{ X_a \}_{a \in B}$
given by the formulae $\psi(x;y)$ and
$\theta(y)$, there is a formula $\phi(z;y)$ so that for each $a \in B$,
the set $\ml L (X_a)$ is defined by $\phi(z;a)$.

\vskip2pt
We recall the following facts about definability in algebraically closed fields.

\begin{fact} \label{factdef} \cite[Theorem A.7]{FLS}  Relative to the theory of algebraically
closed fields (ACF), we have the the following statements.
\begin{enumerate}
\item \label{clodef} The Zariski closure is definable in families.
\item \label{dimdegdef} The dimension and degree of the Zariski
closure of a set are definable in families.
\item \label{irrdef} Irreducibility of the Zariski closure is a definable property.
\item \label{hypdegvardef}
If the Zariski closure is an irreducible hypersurface given by the vanishing of some nonzero polynomial, then the degree of that
polynomial in any particular variable is definable in families.
\item \cite[Lemma 3.5]{FLS}The set of irreducible varieties in $\mathbb A^n$ of dimension $d$ and degree $g$ is a definable family.
\end{enumerate}
\end{fact}

We also need to generalize results on prolongation admissible varieties \cite{FLS} to the partial differential case.
Notations $\tau_l, \nabla_l, B_l$ should be specified beforehand.
For an algebraic variety $X=\ml V(f_1,\ldots,f_o)\subseteq\mb A^n$ defined by polynomials $f_i\in\ml F[y_1,\ldots,y_n]$,
$\tau_l(X)\subseteq\mb A^{n{l+m\choose m}}$ denotes the algebraic variety defined by $(\theta(f_i))_{\theta\in\Theta_{\leq l}}$
considered as algebraic polynomials in $\ml F[\Theta_{\leq l}(\Y)]$  with $\Y=(y_1,\ldots,y_n).$
Thus, $\tau_l\mb A^n=\mb A^{n{l+m\choose m}}$ with coordinates corresponding to variables $(\Y,\Theta_1(\Y),\ldots,\Theta_l(\Y))$.
Given a point $\bar{a}\in\mb A^n$, $\nabla_l(\bar{a})$ denotes the point
$\big(\bar{a},\Theta_1(\bar{a}),\ldots,\Theta_l(\bar{a})\big)\in\tau_l\mb A^n$,
and for a $\Delta$-variety $W\subseteq\mb A^n$,
$B_l(W)$ is the Zariski closure of the set $\{\nabla_l(\bar{a}):\,\bar{a}\in W\}$ in $\tau_l\mb A^n$.
In other words, $B_l(W)=\ml V\big(\mb I(W)\cap\ff[\Theta_{\leq l}(\Y)]\big)\subseteq\tau_l\mb A^n$.

\begin{definition}
Let $V\subseteq\tau_s\mb A^n$ be an algebraic variety over $\ff$.
We say $V$ is {\em prolongation admissible} if  $B_s\big(\mb V(\ml I(V))\big)=V$.
\end{definition}

Irreducible prolongation admissible varieties are of special interest in this paper. 
The following lemma shows that algebraic characteristic sets of irreducible prolongation admissible varieties have a special form,
which generalizes \cite[Lemma 2.13]{FLS} to the partial differential case.
\begin{lemma} \label{lm-prolongationcharset}
Let $V\subseteq\tau_s\mb A^n$ be an irreducible prolongation admissible variety over $\ff$ and
$\mathcal{A}$ a characteristic set of $V$ w.r.t. an ordering induced by some orderly ranking $\mathscr R$ on $\Theta(\Y)$.
For each $k=1,\ldots,n$, let
$$E_k=\big\{\theta y_k\in\lead(\ml A):\,\, (\forall\, \tau y_k\in\lead(\ml A)\wedge (\tau|\theta))\Rightarrow \tau(y_k)=\theta(y_k)\big\}.$$
If $E_k\neq\emptyset$, then for each $\tau y_k\in\Theta_{\leq s}(y_k)$ which is a proper derivative of  some element of $E_k$,
there exists  $A_{\tau,k}\in\ml A$ such that $\lead(A_{\tau,k})=\tau y_k$ and $A_{\tau,k}$ is linear in $\tau y_k$.
\end{lemma}
\proof Follow the proof of \cite[Lemma 2.13]{FLS}. 
Let $W=\mb V \big(\ml I(V)\big)\subseteq\mb A^n$ and $W=\cup_{i=1}^lW_i$ be the irreducible decomposition of $W$.
Since $V$ is prolongation admissible, $B_s(W)=V$.
So there exists some $i_0$ such that $B_s(W_{i_0})=V$.
Suppose $\mathcal B$ is a $\Delta$-characteristic set of $W_{i_0}$ w.r.t. $\mathscr R$.
Let $\ml C=\Theta(\ml B)\cap\ml F[\Theta_{\leq s}(\Y)]$,
$\mathcal C$ is a characteristic set of $B_s(W_{i_0})=V$.
Since $\mathcal C$ and $\ml A$ have the same rank, $\ml A$ should satisfy the desired property.
\qed

\begin{example}
Let $m=2$, $n=1$ and $s=1$. Let $V_1=\mathcal V(y,\delta_1 y, \delta_2y)\subset\tau_1\mathbb A^1$
and $V_2=\mathcal V(y,\delta_2y)\subset\tau_1\mathbb A^1$. Then $\mathbb V(\mathcal I(V_i))=\{0\}\subset\mathbb A^1$ for $i=1,2$.
So $V_1$ is prolongation admissible while  $V_2$ is not prolongation admissible, for $B_1(\mathbb V(\mathcal I(V_2)))=\{(0,0,0)\}\subsetneq V_2=\{(0,c,0): c\in\mathcal E\}$.
This fact can also be obtained by Lemma \ref{lm-prolongationcharset} for $E=\{y\}$ and there  does not exist such an $A_{\delta_1}$ linear in $\delta_1 y$ in a characteristic set of $V_2$.

\end{example}

We now  show that prolongation admissibility is a definable property similarly as in the ordinary differential case \cite[Lemma 2.28]{FLS}.
\begin{lemma} \label{prolongadmissible-definable}
 Let $(V_{{b}})_{{b}\in B}$ be a definable family of algebraic varieties in $\tau_s\mathbb{A}^n$.
 Then $\{b\in B: V_{{b}}\,\text{ is prolongation admissible}\}$ is a definable set.
\end{lemma}
\proof  Follow the proof of \cite[Lemma 2.28]{FLS}.
 Suppose each $V_{b}\subset\tau_s\mathbb{A}^n$  in the definable family  $(V_{{b}})_{{b}\in B}$  is
 defined by $f_i\big({b},(\theta y_j)_{\theta\in\Theta_{\leq s}, 1\leq j\leq n}\big)=0, \,i=1,\ldots,\ell$.
By abuse of notation, let $B_s(V_b)$  be the Zariski closure of $\{\nabla_s(\bar{a}): \nabla_s(\bar{a})\in V_b\}$ in $\tau_s\mathbb{A}^n$.
Then $\deg(B_s(V_b))$ has a uniform bound $T$ in terms of the degree bound $D$ of the $f_i$, $m$, $n$, $\ell$ and $s$.
Indeed, let $z_{j,\theta}\,(j=1,\ldots,n;\theta\in\Theta_{\leq s})$ be new $\Delta$-variables and replace $\theta(x_j)$ by $z_{j,\theta}$ in each $f_i$
to get a new differential polynomial   $g_{i}$.
 Consider the new differential  system $S:=\{ g_1,\ldots,g_\ell, \delta_k(z_{j,\theta})-z_{j,\delta_k\theta}:\,\,k=1,\ldots,m;\theta\in\Theta_{\leq s-1}\}$.
 Regard $S$ as a pure algebraic polynomial system in $z_{j\theta}$ and $\delta_k(z_{j\theta})$ temporarily, and let $U$ be the
 Zariski closed set  defined by $S$ in $\tau_1(\tau_{s}\mathbb{A}^n)$.
 Let $Z=\{\bar{c}\in\tau_s\mathbb{A}^n:  \nabla(\bar{c})\in U\}$.
  Clearly,  $Z=\{\nabla_s(\bar{a}): \nabla_s(\bar{a})\in V_b\}$.
By \cite[Corollary 4.5]{FS}, the degree of the Zariski closure of $Z$, namely $B_s(V_b)$, is bounded by  some number $D_1$ which depend on $D$, $m$, $n$, $\ell$ and $s$.

By \cite[Proposition 3]{heintz}, an irreducible algebraic variety $V\subseteq\tau_s\mb A^n$ can be defined by $n{s+m\choose m}+1$ polynomials of degree bounded by the degree of $V$. Since the degree of an algebraic variety is the sum of the degrees of its components, combined with Kronecker's theorem \cite[p.146]{ritt},
$B_s(V_b)$ could be defined by at most $n{s+m\choose m}+1$ polynomials of degree bounded by $D_1$.
Hence, $(B_s(V_b))_{b\in B}$ is a definable family.
Since $V_b$ is prolongation admissible if and only if $V_b=B_s(V_b)$,
which implies that $\{b: V_b \,\text{is  prolongation admissible}\}$ is a definable set.
\qed

\begin{definition}
Let $V\subset\tau_l\mb A^n$ be an irreducible prolongation admissible  variety
and $W=\mb V(I(V))$ be the $\Delta$-variety defined by defining equations of $V$.
A component $W_1$ of $W$ is called a {\em dominant component} if $B_s(W_1)=V$.
\end{definition}

The following result shows how to get the desired unique irreducible $\Delta$-varieties from irreducible prolongation admissible varieties,
 where additional conditions are required to generalizes \cite[Lemma 2.15]{FLS} to the partial differential case.

\begin{lemma} \label{lem-uniquedominant}
Let $V\subseteq\tau_s(\mb A^{n})$ be an irreducible prolongation admissible variety of dimension $(d+1){s+m\choose m}-1$
and in the case  $s>0$, suppose additionally $\pi_{s,0}(V)$ is of dimension $d+1$.
Then $W=\mb V(I(V))$ has a unique dominant component $W_1$ and $\omega_{W_1}(t)=(d+1){t+m\choose m}-{t+m-s\choose m}$.
\end{lemma}
\proof Two cases should be considered according to whether $s=0$ or not.

Case 1) $s=0$. In this case, $\ml P=I(V)$ is a prime ideal of $\ff[y_1,\ldots,y_n]$ of dimension $d$.
By \cite[p.200, Proposition 10]{kolchin}, $\{\ml P\}$ is a prime $\Delta$-ideal of $\ff\{y_1,\ldots,y_n\}$ with
Kolchin  polynomial $\omega_{\{\ml P\}}(t)=d{t+m\choose m}$.
Thus, $W=\mb V(\ml P)$ itself is its dominant component  and satisfies the desired property.

Case 2) $s>0$.
Fix an orderly ranking $\mathscr R$ on $\ff\{\Y\}$ and denote $\mathscr R_l$ to be the ordering on
$\Theta_{\leq l}(\Y)$ induced by $\mathscr R.$
Since $\pi_{s,0}(V)$ is of dimension $d+1$, a characteristic set of  the Zariski closure of ${\pi_{s,0}(V)}$ w.r.t. $\mathscr R_0$ is of the form
$B_1,\ldots,B_{n-d-1}$ where $\lead(B_i)=y_{\sigma_i}$ for each $i$.
Since $V$ is irreducible and prolongation admissible, by Lemma \ref{lm-prolongationcharset},
 $S=\{\theta(y_{\sigma_i}): \,\,\ord(\theta)\leq s, i=1,\ldots,n-d-1\}$ is a subset of the leaders of
 a characteristic set $\mathcal A$ of $V$ w.r.t. $\mathscr R_s$.
Since the dimension of $V$ is  $(d+1){s+m\choose m}-1$,  $\lead(\mathcal A)=S\cup\{\tau(y_{\sigma_{n-d}})\}$
for some $\tau\in\Theta_s$ and $\sigma_{n-d}\in\{1,\ldots,n\}\backslash\{\sigma_1,\ldots,\sigma_{n-d-1}\}$.
So there exists $B_{n-d}\in\mathcal A$ s.t. $\lead(B_{n-d})=\tau(y_{\sigma_{n-d}})$.

Let $\ml B =<B_1,\ldots,B_{n-d}>$.
Clearly, $\ml B$ is an irreducible coherent  autoreduced set of $\ff\{y_1,\ldots,y_n\}$, by \cite[Lemma 2, p.167]{kolchin},
$\ml B$ is a $\Delta$-characteristic set of a prime $\Delta$-ideal $\ml P\subseteq \ff\{y_1,\ldots,y_n\}$ w.r.t. $\mathscr R$.
Clearly, $\ml P=\sat(\ml B)$ and its Kolchin  polynomial $\omega_{\ml P}(t)=(d+1){t+m\choose m}-{t+m-s\choose m}$.
We now show that  $\mb V(\ml P)\subseteq W$ and $B_s(\mb V(\ml P))=V$.
Since $V$ is an irreducible prolongation admissible variety, there exists a point $\bar{a}\in\mb A^n$
such that $\nabla_{s}(\bar{a})$ is a generic point of $V$.
So as $\Delta$-polynomials, $B_{i}$ vanishes at $\bar{a}$ while $H_{\ml B}$ does not.
Thus, $\ml P$ vanishes at $\bar{a}$, and consequently, $\nabla_{s}(\bar{a})\in B_s(\mb V(\ml P))$.
So $V\subseteq B_s(\mb V(\ml P))$.
Since both $V$ and $B_s(\mb V(\ml P))$ are irreducible varieties of the same dimension, $B_s(\mb V(\ml P))=V.$
So, $\ml I(V)=\ml P\cap\ff[\Theta_{\leq s}(\Y)]\subseteq\ml P$, as a consequence, $\mb V(\ml P)\subseteq \mb V(I(V))=W$.

Suppose $W_0$ is a dominant component of $W$.
Given a generic point $\xi\in W_0$, $\nabla_s(\xi)$ is a generic point of $V$.
So, $\ml B$ vanishes at $\xi$ and $H_{\ml B}$ does not vanish at $\xi$.
Thus, $\mb V(\ml P)$ vanishes at $\xi$, i.e., $W_0\subseteq \mb V(\ml P)$.
So  $W_0=\mb V(\ml P)$. Thus,  $\mb V(\ml P)$ is the unique dominant component $W$ and $\omega_{\mb V(\ml P)}(t)=(d+1){t+m\choose m}-{t+m-s\choose m}$.
\qed

\subsection{Proof of the main theorem}
Before proving the main theorem, we need to bound the degree of $B_s(V)$ similar as in \cite[Proposition 4.3]{FLS} to get the candidates of the algebraic Chow varieties
which can be used to paramertrize $\Delta$-varieties in $G_{(n,d,s,r)}(\ff)$.

\begin{lemma} \label{lem-bound}
Let $V\subset\mb A^n$ be an irreducible $\Delta$-variety in $G_{(n,d,s,r)}(\ff)$.
Then $B_s(V)$ is an irreducible  variety in $\tau_s(\mb A^n)$ over $\ff$ of dimension $(d+1){s+m\choose m}-1$
and the degree of $B_s(V)$ satisfies
  $$r\Big/{s+m\choose m}\leq \deg(B_s(V))\leq \big[(s+1)(d+1)r\big]^{n(s+1){s+m\choose m}+1}.$$
\end{lemma}
\proof It is clear that $B_s(V)$ is an irreducible  variety in $\tau_s(\mb A^n)$ of dimension $(d+1){s+m\choose m}-1$.
For the degree bound,   we will first show  that $ \deg(B_s(V))\leq \big[(s+1)(d+1)r\big]^{n(s+1){s+m\choose m}+1}$.
Since $V\in G_{(n,d,s,r)}(\ff)$,  the $\Delta$-Chow form $F(\bu_0,\ldots,\bu_d)$ of $V$ exists,
and $\ord(F)=s$, $\deg(F,\bu_0^{[s]})=r$.
Let   $\ml J=[\mb I(V),\mb  L_0,$ $\ldots,\mb L_d]\subseteq{\ff\{\bu_0,\ldots,\bu_d, \Y\}}$.
Let $\mr R$ be a ranking  on $\ff\{\bu_0,\ldots,\bu_d, \Y\}$ satisfying 1) $\theta(u_{ij})<\tau(y_k)$ for any $\theta$ and $\tau$,
 and 2) $\mathscr R$ restricted to $\bu_0,\ldots,\bu_d$ is an orderly ranking with $u_{00}$ greater than any other $u_{ij}$.
Let $\mr R_s$ be the ordering on $\bu_0^{[2s]},\ldots,\bu_d^{[2s]}$ and $\Y^{[s]}$  induced by $\mr R$.
Suppose $\lead(F)=\theta(u_{00})$ for some $\theta\in\Theta_s$ and $\sept_F=\frac{\partial F}{\partial \theta(u_{00})}$.

 By Theorem \ref{th-relationchow}, for $j=1,\ldots,n$,
  the polynomial $G_{j}=\sept_Fy_j-\frac{\partial F}{\partial \theta(u_{0j})}\in \ml J$ and note that $\deg(G_j)=(d+1)r$.
We construct polynomials $G_{j,\theta}\in \ml J$ for $\theta\in\Theta_{\leq s}$ with $\rank(G_{j,\theta})=\theta(y_j)$ and $\deg(G_{j,\theta})\leq \big({\ord(\theta)}+1\big)(d+1)r$ inductively on the order of $\theta$.
 Set $G_{j,1}=G_j$.
 Let $G_{j,\delta_i}=\rem(\delta_i(G_{j,1}),G_{j,1})$ be the algebraic remainder of $\delta_i(G_{j,1})$ with respect to $G_{j,1}$.
Clearly, $G_{j,\delta_i}\in \ml J$ and is of the form $G_{j,\delta_i}=\sept_F^2\delta_i(y_j)+T_{j,\delta_i}$ for some  $T_{j,\delta_i}\in\ff[\bu^{[s+1]}]$.
An easy calculation shows that $\deg(G_{j,\delta_i})\leq 2(d+1)r$.
Suppose the desired $G_{j,\tau}=\sept_F^{\ord(\tau)+1}\tau(y_j)+T_{j,\tau}\,(\tau\in\Theta_{\leq k})$ have been constructed,
we now define $G_{j,\tau}\,(\tau\in\Theta_{k+1})$.
 For $\tau\in\Theta_{k+1}$,
 let $G_{j,\tau}$ be the algebraic remainder of $\tau(G_j)$ with respect to $<G_{j,\tau}:\,\tau\in\Theta_{\leq k}; k\leq s>$.
 Then $G_{j,\tau}\in  \ml J$ and $G_{j,\tau}=\sept_F^{k+2}\tau(y_j)+T_{j,\tau}$ ,where $T_{j,\tau}\in\ff[\bu^{[k+1+s]}]$
 satisfies $\deg(T_{j,\tau})\leq (k+2)(d+1)r$.
 In  this way, the polynomials $G_{j,\tau}\in \ml J\,(\tau\in\Theta_{\leq s})$ are constructed.

Clearly,  $<F^{[s]},  G_{j,\tau}:\,\tau\in \Theta_{\leq s}>$ is an irreducible ascending chain under $\mr R_s$,
so $\ml J_s=\big(F^{[s]}, (G_{j,\tau})_{\tau\in \Theta_{\leq s}}\big):\sept_F^\infty$ is a prime ideal
in $\ff[\Y^{[s]},\bu_0^{[2s]}, \ldots,\bu_d^{[2s]}]$, which is a component of $\ml V(F^{[s]},(G_{j,\tau})_{\tau\in \Theta_{\leq s}})$.
By Be\'{z}out Theorem \cite[Theorem 1]{heintz},  we have
\begin{eqnarray}
\deg(\ml J_s)&\leq& [(d+1)r]^{{s+m\choose m}}\cdot\prod_{j=1}^n \prod_{\theta\in\Theta_{\leq s}}\deg(G_{j\theta})\nonumber\\
&\leq& [(d+1)r]^{{s+m\choose m}}\cdot\prod_{j=1}^n\prod_{l=0}^s[(l+1)(d+1)r]^{{s+m\choose m}}\nonumber\\
&<& [(s+1)(d+1)r]^{n(s+1){l+m\choose m}+1}. \nonumber
\end{eqnarray}
Let $\ml J_s'=\ml J_s\cap\ff[\Y^{[s]}]$.
We claim that $\ml J_s'=\mb I(V)\cap\ff[\Y^{[s]}]$.
Indeed, on the one hand, $\ml J_s'\subseteq \ml J\cap\ff[\Y^{[s]}]=\mb I(V)\cap\ff[\Y^{[s]}]$;
on the other hand, for any polynomial $H\in\mb I(V)\cap\ff[\Y^{[s]}]$,
the algebraic remainder of $H$ with respect to $\langle G_{j,\tau}:\,\tau\in \Theta_{\leq s}\rangle$
is a polynomial $H_1\in \ml J\cap\ff\{\bu_0,\ldots,\bu_d\}=\sat(F)$ with $\ord(H_1)\leq 2s$.
Thus, $H_1\in \asat(F^{[s]})$ and $H\in\ml J_s$.
So  by \cite[Lemma 2]{heintz} or \cite[Theorem 2.1]{issac2011}, $\deg(\mb I(V)\cap\ff[\Y^{[s]}])=\deg(\ml J_s')\leq \deg(\ml J_s)$.

Now, we show $\deg(\mb I(V)\cap\ff[\Y^{[s]}])\geq r/{s+m\choose m}$.
By item 3) of Theorem \ref{th-relationchow},  $\big[\ml I\cap\ff[\Y^{[s]}],$ $\mb  L_0^{[s]},\ldots,\mb L_d^{[s]}\big]\cap\ff[\bu_0^{[s]},\ldots,\bu_d^{[s]}]=(F)$. Similar to the procedures in \cite[Theorem 6.25]{FOCM}, the $\Delta$-Chow form of $\mb I(V)$ could be obtained from the algebraic Chow form of $\mb I(V)\cap\ff[\Y^{[s]}]$
by algebraic specializations.
So $(d+1)r\leq  (d+1){s+m\choose m}\deg(\mb I(V)\cap\ff[\Y^{[s]}])$ and $\deg(B_s(V))=\deg(\mb I(V)\cap\ff[\Y^{[s]}])\geq r/{s+m\choose m}$.
\qed

\begin{remark}
 In the ordinary differential case \cite[Proposition 4.3]{FLS}, the construction of $G_{j,k}$ is much easier
 and each $G_{j,k}\,(k\leq s)$ could be chosen from $\ff[\bu_0^{[s]},\ldots,\bu_d^{[s]},\Y^{[s]}]$.
 However, due to the more complicated structure of partial differential characteristic sets, 
 the method could not be adapted here because there may exist $\tau\in\Theta_{\leq s}$ 
 such that any derivative of $\tau(u_{00})$ does not appear in $F$.
 Also,  here $G_{j,\theta}\in\ff[\bu_0^{[2s]},\ldots,\bu_d^{[2s]},\Y^{[s]}]$ for $\theta\in\Theta_{\leq s}$.
\end{remark}

Now, we are ready to prove that $\Delta$-Chow varieties of $\mb A^n$ of  index $(d,s,r)$  exist for all $n, d, s, r$.
As mentioned in the beginning of this section, we will use certain algebraic Chow varieties to parametrize $\Delta$-varieties in $G_{(n,d,s,r)}$.
For the sake of later use, we shall briefly recall the concept of algebraic Chow varieties here.

For an irreducible variety $V\subseteq\mb P^n$ over $\ff$ of dimension $d$, the algebraic Chow form of $V$ is the polynomial $H(\bu_0,\ldots,\bu_d)$ whose vanishing gives a necessary and sufficient condition for $V$ and $d+1$ hyperplanes $u_{i0}y_0+\sum_{j=1}^nu_{ij}y_j=0$ having a nonempty intersection in $\mb P^n$. Here $\bu_i=(u_{i0},\ldots,u_{in})$. For a $d$-cycle $W$ in $\mb P^n$, $W=\sum_{i=1}^lt_iW_i$ with $t_i\in\mb N$ and $\dim(W_i)=d$, the Chow form of $W$ is the product of Chow forms of $W_i$ with multiplicity $t_i$.
Its degree in each $\bu_i$ is called the degree of $W$ and its coefficient vector, regarded as a point in a higher dimensional projective space,  is defined to be the {\em Chow coordinates} of $W$. 
The  set of Chow coordinates of all $d$-cycles in  $\mathbb P^n$ of  degree $e$ is a projective variety in the Chow coordinate space \cite{chow, hodge}, called the Chow variety (of index $(d,e)$). 
Moreover, by \cite[p. 57, Theorem II]{hodge} and \cite[p. 697-698 ]{chow}, the defining equations of  Chow varieties are homogenous polynomial equations over $\mathbb Q$.
However, the affine Chow variety  of all $d$-cycles in $\mb A^n$ of degree $e$ is not Zariski closed in the Chow coordinate space,  but it is always a constructible set \cite[Proposition 3.4]{FLS}, denoted  by $\text{Chow}_{n}(d,e)$. 
Each point of $\text{Chow}_{n}(d,e)(\ff)$ represents a $d$-cycle in $\mb A^n$ of degree $e$ over $\ff$. 
All the Chow varieties we use here are affine ones.

To show the existence of $\Delta$-Chow varieties of $\mb A^n$ of  index $(d,s,r)$, two cases should be considered separately.
In the case $s=0$,  the $\Delta$-Chow form of each $V\in G_{(n,d,0,r)}$ is just equal to the
Chow form of  $B_0(V)\subseteq \mb A^n$, so the set of $\Delta$-Chow coordinates of $\Delta$-varieties in $G_{(n,d,0,r)}$ is just the same as the set of Chow coordinates of all irreducible varieties in $\mathbb A^n$ of dimension $d$ and degree $r$. By item 5) of Fact \ref{factdef}, the latter set is a definable subset of  $\text{Chow}_{n}(d,r)$, so $G_{(n,d,0,r)}$ is a constructible set.
Below, we focus on the case  $s>0$.

Let $\text{Chow}_{n{s+m\choose m}}\big((d+1){s+m\choose m}-1, e\big)$ be the affine Chow variety
of all cycles  in $\tau_{s}(\mb A^n)$  of dimension $(d+1){s+m\choose m}-1$ and degree $e$.
Consider the disjoint union of algebraic constructible sets $$\mathcal{C}=\bigcup_{D_1\leq e \leq D_2}{\chow}_{n{s+m\choose m}}\big((d+1){s+m\choose m}-1,e\big)$$
where $D_1, D_2$ are the lower and upper bounds given in Lemma \ref{lem-bound}.
So each point $a\in \mathcal{C}$ represents a $[(d+1){s+m\choose m}-1]$-cycle in $\tau_s\mathbb{A}^n$.
To represent an irreducible $\Delta$-variety $V$ of the desired Kolchin  polynomial and $\Delta$-degree by a point in $\mathcal{C}$,
we only need to consider those irreducible varieties with Chow coordinates in $\mathcal{C}$.

Let $\mathcal{C}_1$ be the subset consisting of all points $a \in\mathcal{C}$ such that
 $a$  is the Chow coordinate of an irreducible variety $W$ which  is  prolongation admissible and additionally satisfies the following conditions:
\begin{enumerate}
\item[1)] $\pi_{s,0}(W)$ is of dimension $d+1$;
\vskip5pt
\item[2)] the unique  dominant component of the $\Delta$-variety defined by equations of $W$ is of $\Delta$-degree $g$.
\end{enumerate}

\begin{theorem}\label{deltachowexists}
The set $\ml C_1$ is a $\Delta$-constructible set and
the map which associates an irreducible $\Delta$-variety $V\subseteq\mb A^n$ in $G_{(n,d,s,r)}$ with the Chow coordinate of
the irreducible  variety $B_s(V)\subseteq\tau_l(\mb A^n)$ identifies $G_{(n,d,s,r)}$ with
$\ml C_1$.  In particular, the $\Delta$-Chow variety of all irreducible $\Delta$-varieties of Kolchin  polynomial $(d+1){t+m\choose m}-{t+m-s\choose m}$  and $\Delta$-degree $r$ exists.
\end{theorem}

\proof
First, we show $\ml C_1$ is a $\Delta$-constructible set.
From the definition of Chow coordinates, we know each ${\chow}_{n{s+m\choose m}}\big((d+1){s+m\choose m}-1,e\big)$
actually represents a definable family  $S_e:=(F_c)_{c\in  {\chow}_{n{s+m\choose m}}\big((d+1){s+m\choose m}-1,e\big)}$ of
homogenous polynomials which are Chow forms of algebraic cycles in $\tau_s\mathbb A^{n}$ of dimension
$(d+1){s+m\choose m}-1$  and degree $e$.
The algebraic cycle whose  Chow coordinate is $c$  is irreducible if and only if its Chow form $F_c$ is irreducible.
Since irreducibility is a definable property, the set $\ml C_0=\big\{c\in  {\chow}_{n{s+m\choose m}}((d+1){s+m\choose m}-1,e): F_c\,\text{is irreducible}\big\}$ is a definable set.
 Take an arbitrary $c\in \ml C_0$ and the corresponding polynomial $F_c\in S_e$ for an example.
Let $V_c$ be the corresponding irreducible variety with Chow coordinate $c$.
By item 5) of Fact \ref{factdef}, $(V_{c})_{c\in \ml C_0}$ is a definable family (For details, please refer to the proof of \cite[Lemma 3.5]{FLS}).
And by Lemma \ref{prolongadmissible-definable} and Fact \ref{factdef},
$\ml C_2=\big\{c\in \ml C_0: V_c \,\text{ is prolongation admissible and}\, \dim(\pi_{s,0}(V_c))=d+1\big\}$ is a definable set.
Then by  Lemma \ref{lem-uniquedominant}, for each $c\in\ml C_2$, the $\Delta$-variety corresponding to $V_c$  has a unique dominant component $W_c$ and the Kolchin  polynomial of $W_c$ is $(d+1){t+m\choose m}-{t+m-s\choose m}$.

Since the Kolchin  polynomial of $W_c$ is $(d+1){t+m\choose m}-{t+m-s\choose m}$,
the $\Delta$-Chow form of $W_c$ exists.
Let $U$ be the algebraic variety in $\tau_s\mathbb{A}^n\times \big(\mathbb{P}^{(n+1){s+m\choose m}-1}\big)^{d+1}$ defined by the defining formulae of $V_c$ and
$\theta(L_i)=0$ for $\theta\in\Theta_{\leq s}$ and $i=0,\ldots,d$ with each $\theta(L_i)=\theta(u_{i0})+\sum_{k=1}^n\sum_{\tau|\theta}{\theta\choose \tau}\frac{\theta}{\tau}(u_{ik})\tau(y_k)$ regarded as a polynomial in  variables $\Theta_{\leq s}(y_{k})$ and $\Theta_{\leq s}(u_{ik})$.
Since $B_s(W_{c})=V_{c}$, by item 3) of Theorem \ref{th-relationchow},
the Zariski closure of the image of $U$ under the following projection map
 $$\pi:\,\tau_s\mathbb{A}^n\times \big(\mathbb{P}^{(n+1){s+m\choose m}-1}\big)^{d+1}\longrightarrow \big(\mathbb{P}^{(n+1){s+m\choose m}-1}\big)^{d+1}$$ is an irreducible variety of codimension 1,
and the defining polynomial $F$ of $\overline{\pi(U)}$ is the $\Delta$-Chow form of $W_{c}$.
By item 4) of Fact \ref{factdef}, the total degree of $F$ is definable in families; this quantity is just equal to $(d+1)$ times the $\Delta$-degree  of  $W_{c}$.
So  the $\Delta$-degree of $W_{c}$ is definable in families.
Hence, $\mathcal{C}_1$ is a definable set, and also a $\Delta$-constructible set due to the fact that the theory
DCF$_{0,m}$ eliminates quantifiers \cite[Theorem 3.1.7]{mcgrail}.

By Lemma~\ref{lem-uniquedominant} and its proof,
each irreducible variety $V$ corresponding to a point of $\ml C_1$ determines an irreducible  $\Delta$-variety
$W\in G_{(n,d,s,r)}$, where $W$ is the unique dominant component of the $\Delta$-variety corresponding to the prolongation admissible variety $V$.
And on the other hand, each $W\in G_{(n,d,s,r)}$ determines the corresponding algebraic irreducible variety $B_s(W)$,
whose Chow coordinate is  a point of  $\ml C_1$ guaranteed by Lemma~\ref{lem-bound}.
So we have established a natural one-to-one correspondence between $G_{(n,d,s,r)}$ and ${\mathcal C}_1$.
Thus, $G_{(n,d,s,r)}$ is represented by the $\Delta$-constructible set  $ {\mathcal C}_1$. \qed

  \section{ Conclusion}
In this paper,  a quasi-generic partial differential intersection theorem is first given. Namely,  the intersection of an irreducible partial differential variety $V$ with a  quasi-generic differential hypersurface of order $s$ is shown to be an irreducible differential variety with  Kolchin   polynomial $\omega_V(t)-{t-s+m\choose m}$.
Then partial differential Chow forms are defined for irreducible partial differential varieties of Kolchin  polynomial $(d+1){t+m\choose m}-{t+m-s\choose m}$
and basic properties similar to their algebraic and ordinary differential counterparts are presented.
Finally, differential Chow coordinate representations are defined for  such partial differential varieties,
and the set of all irreducible partial differential varieties of fixed Kolchin  polynomial and differential degree
is shown to have a structure of differentially constructible set.

The above results have generalized  the theory of differential Chow forms and Chow varieties obtained for the ordinary differential case \cite{GLY, FLS}  to their partial differential analogs.
However, the theory of partial differential Chow forms and partial differential Chow varieties is far from well-developed and there are several unsolved problems.
As stated in Conjecture \ref{conjecture}, we conjecture that $\omega_V(t)=(d+1){t+m\choose m}-{t+m-s\choose m}$ for some $d,s\in\mb N$ is not only a sufficient condition,
but also a necessary condition such that the partial differential Chow form of  $V$ exists.
Another problem is how to represent general (irreducible) partial differential varieties by coordinates and further how to provide a set of partial differential varieties of fixed characteristics
with a structure of partial differential constructible set.

 \section{Acknowledgements}
The author is grateful to the anonymous referee for the helpful comments and constructive suggestions on a previous version of this manuscript.
This work was supported by NSFC under Grants (11971029, 11688101, 11671014, 11301519).


\begin{thebibliography}{10}


\bibitem{brownawell}
Brownawell, W.~D. (1987).
Bounds for the degrees in the Nullstellensatz.
{\em Ann. Math.} {126(3)}:577--591.
  DOI: {10.2307/1971361}.
  
  
\bibitem{BC}
Buium, A. and Cassidy, P. (1998).
Differential algebraic geometry and differential  algebraic groups.
{\em In H. Bass et al eds, Selected Works of Ellis Kolchin, with Commentary,
  American Mathematical Society, Providence, RI}, pp. 567--636.

\bibitem{chow}
Chow, W.-L. and van~derWaerden, B.~L. (1937). Zur algebraischen {G}eometrie. {IX}.
{\em Math. Ann.} { 113(1)}: 692--704.
DOI: {10.1007/BF01571660}.


\bibitem{Eisenbud-Harris}
Eisenbud, D. and Harris, J. (2016).
3264 and all that---a second course in algebraic geometry.
Cambridge University Press, Cambridge.

\bibitem{eisenbud}
Eisenbud, D., Schreyer, F.-O., and Weyman, J. (2003).
Resultants and {C}how forms via  exterior syzygies.
{\em J. Amer. Math. Soc.} {16(3)}: 537--579.
DOI:{10.1090/S0894-0347-03-00423-5}.

\bibitem{freitag}
Freitag, J. (2015).
Bertini theorems for differential algebraic geometry.
{\em ArXiv:1211.0972v3}.


\bibitem{FS}
Freitag, J. and Le\'{o}n~S\'{a}nchez, O. (2016).
Effective uniform bounding in partial differential fields.
{\em Adv. Math.} {288}:308--336.
DOI:{10.1016/j.aim.2015.10.013}.

\bibitem{FLS}
Freitag, J., Li, W., and Scanlon, T. (2017).
Differential Chow varieties exist. With an appendix by William Johnson.
{\em J. Lond. Math. Soc. (2)} { 95(1)}:128--156.
DOI: {10.1112/jlms.12002}.

\bibitem{GLY}
Gao, X.-S., Li, W., and Yuan, C.-M. (2013).
Intersection theory in differential algebraic geometry: generic intersections and the differential {C}how form.
{\em Trans. Amer. Math. Soc.} {365(9)}:4575--4632.
DOI: {10.1090/S0002-9947-2013-05633-4}.

\bibitem{harris1992algebraic}
Harris, J. (1992). Algebraic {G}eometry: {A} {F}irst {C}ourse, GTM {133}. Springer, New York.

\bibitem{heintz}
Heintz, J. (1983).
Definability and fast quantifier elimination in algebraically  closed fields.
{\em Theoret. Comput. Sci.} { 24(3)}:239--277.
DOI: {10.1016/0304-3975(83)90002-6}.

\bibitem{hodgevol1}
Hodge, W. and Pedoe, D. (1994).
Methods of algebraic geometry. {V}ol. {I}.
Cambridge University Press, Cambridge.

\bibitem{hodge}
Hodge, W. V.~D. and Pedoe, D. (1994).
Methods of algebraic geometry, {V}ol. {II}.
Cambridge University Press, Cambridge.

\bibitem{Complexitychowform}
Jeronimo, G., Krick, T., Sabia, J., and Sombra, M. (2004).
The computational complexity of the {C}how form.
{\em Found. Comput. Math.} { 4(1)}: 41--117.
DOI: {10.1007/s10208-002-0078-2}.

\bibitem{Kolchin1964}
Kolchin, E.~R. (1964).
The notion of dimension in the theory of algebraic
  differential equations.
{\em Bull. Amer. Math. Soc.} {70}:570--573.
  DOI: {10.1090/S0002-9904-1964-11203-9}.

\bibitem{kolchin}
Kolchin, E.~R. (1973).
Differential algebra and algebraic groups.
Academic Press, New York-London.

\bibitem{KLMP}
Kondratieva, M.~V., Levin, A.~B., Mikhalev, A.~V., and Pankratiev, E.~V. (1999).
Differential and difference dimension polynomials.
 Mathematics and its Applications, volume {461}.
Kluwer Academic Publishers, Dordrecht.

\bibitem{sanchez}
Le\'{o}n~S\'{a}nchez, O. (2012).
Geometric axioms for differentially closed fields with several commuting derivations.
{\em J. Algebra} {362}, 107--116.
DOI: {10.1016/j.jalgebra.2012.03.043}.

\bibitem{issac2011}
Li, W., Gao, X.-S., and Yuan, C.-M. (2011).
Sparse differential resultant.
In I{SSAC} 2011---{P}roceedings of the 36th {I}nternational {S}ymposium on
  {S}ymbolic and {A}lgebraic {C}omputation, ACM, New York, June 8-11 (San Jose, California, USA), pp. 225--232.

\bibitem{LG}
Li, W. and Gao, X.-S. (2012).
Chow form for projective differential variety.
{\em J. Algebra} { 370}, 344--360.
 DOI: {10.1016/j.jalgebra.2012.07.047}.


\bibitem{FOCM}
Li, W., Yuan, C.-M., and Gao, X.-S. (2015).
Sparse differential resultant for {L}aurent differential polynomials.
{\em Found. Comput. Math.} {15(2)}:451--517.
  DOI: {10.1007/s10208-015-9249-9}.

\bibitem{mcgrail}
McGrail, T. (2000).
The model theory of differential fields with finitely many commuting derivations.
{\em J. Symbolic Logic} { 65(2)}: 885--913.
DOI: {10.2307/2586576}.

\bibitem{philippon}
Philippon, P. (1986).
Crit\`eres pour l'indpendance algbrique.
{\em Inst. Hautes \`Etudes Sci. Publ. Math.} { 64}: 5--52.

\bibitem{pierce}
Pierce, D. (2014).
Fields with several commuting derivations.
{\em J. Symb. Log.} {79(1)}:1--19.
  DOI: {10.1017/jsl.2013.19}.

\bibitem{ritt}
Ritt, J.~F. (1950).
Differential {A}lgebra.
American Mathematical Society, New York, N. Y..

\bibitem{wu}
Wu, W.-T. (1989).
On the Foundation of Algebraic Differential Polynomial Geometry.
{\em Sys. Sci. Math. Sci.} { 2(4)}:{289--312}.


\end{thebibliography}
\end{document}